\documentclass[12pt]{amsart}
\usepackage{amsmath}
\usepackage{amssymb}
\usepackage{graphicx}
\newtheorem{theorem}{Theorem}
\newtheorem{conjecture}[theorem]{Conjecture}
\newtheorem{lemma}[theorem]{Lemma}
\newtheorem{corollary}[theorem]{Corollary}

\numberwithin{equation}{section}
\numberwithin{theorem}{section}

\begin{document}
\title{Autocorrelation of ratios of $L$-functions}

\author {Brian Conrey \and David W.~Farmer \and Martin R. Zirnbauer}

\abstract
We give a new heuristic for all of the main terms in the
quotient of products of $L$-functions averaged over a family.
These conjectures generalize the recent conjectures
for mean values of $L$-functions.  Comparison is made to the
analogous quantities for the characteristic polynomials of
matrices averaged over a classical compact group.
\endabstract

\maketitle

\section{Introduction}

Conjectures for the moments of $L$-functions have been sought for many
decades, with very little progress until the introduction of
random matrix theory into the subject~\cite{CF, KS1, KS2}.
The predictions using random matrix theory provide plausible
conjectures for the \emph{leading order} asymptotics of the
moments of a family of $L$-functions.
More recently, conjectures for all of the main terms of the
moments have been found, using heuristics based on the
harmonic detector of a family~\cite{CFKRS2}, and also using
a plausible conjecture for multiple Dirichlet series~\cite{DGH}.
Although the more refined conjectures do not make use of random matrix
theory, they are supported by the similarity to the analogous
theorems for random matrices~\cite{CFKRS1}.

In this paper we generalize the heuristic method of~\cite{CFKRS2}
to the case of ratios of products of $L$-functions.  In our
companion paper~\cite{CFZ} and in ~\cite{HPZ} the methods of supersymmetry
are used
to evaluate the analogous quantity for the characteristic polynomials
of matrices averaged over the compact classical groups.  We quote
those results in Section~\ref{sec:results} for comparison with
our conjectures given in Section~\ref{sec:conjectures}.

The usefulness of averages of ratios was first suggested by
Farmer~\cite{F1}, who made the following conjecture
about the Riemann zeta-function. For complex numbers
$\alpha, \beta, \gamma, \delta$ with  real parts that are positive
and of size $c/\log T$,
\begin{equation}\label{eqn:farmerconjecture}
\frac{1}{T}\int_0^T \frac{\zeta(s+\alpha)\zeta(1-s+\beta)}
{\zeta(s+\gamma)\zeta(1-s+\delta)}~dt \sim
1+(1-T^{-\alpha-\beta})\frac{(\alpha-\gamma)(\beta-\delta)}{(\alpha+\beta)
(\gamma+\delta)}.
\end{equation}
This conjecture, developed in conjunction
with the conjecture on  `long mollifiers,'
had a somewhat shaky justification. Yet
implied many things of great interest about the Riemann
zeta-function, such as the pair correlation conjecture of Montgomery~\cite{M},
Levinson's formula for a `mollified' mean square
of $\zeta(s)$,
and asymptotic formula for moments of the logarithmic
derivative of the zeta-function near
the critical line.
Moreover, it satisfied many consistency checks
involving exchanging parameters, using the functional equation
of $\zeta(s)$ and letting variables approach 0 or $\infty$.
Farmer later extended the conjecture to include similar
integrals but with ratios of up to three zeta-functions in the
numerator and denominator. These formulas were also found to
imply interesting statements about the zeros of
$\zeta(s)$, including the triple correlation conjecture
of Hejhal~\cite{He} and Rudnick-Sarnak~\cite{RS}.

In another direction, Goldston and Montgomery~\cite{GM} have
shown that the pair correlation conjecture is
equivalent to a statement about the second moment of
primes in short intervals. In~\cite{GGM} an asymptotic formula
for the mean-square of the logarithmic derivative of $\zeta(s)$
was shown to be equivalent to both the pair correlation and
(hence) the second moment of primes.  Conjecture~(\ref{eqn:farmerconjecture})
encapsulates these results.


Thus, it was of great interest that the analogue of
formula~(\ref{eqn:farmerconjecture}) was found to be true
for the characteristic polynomial of a matrix averaged over the
unitary group~$U(N)$.  This was first observed by Nonenmacher and
Zirnbauer at a workshop at MSRI in 1999.
Since it is believed that families of $L$-functions can be modeled by
the characteristic polynomials from one of the classical compact groups,
these formulas for other compact groups lead  to conjectures for the averages,
over a family,
of ratios of products of $L$-functions.
In every case where we can prove and asymptotic formula, or conjecture one
from number theoretic reasoning,
we have agreement with the conjectures presented here.
We noted above that the ratios conjectures imply Montgomery's
pair correlation conjecture.  The ratios conjectures contain much
more information and can be used to make very precise conjectures about
the distribution of zeros of $L$-functions.  In addition to the examples
we give in Section~\ref{sec:applications} of this paper,
Conrey and Snaith~\cite{CS} have given a large number of applications.

In Section~\ref{sec:notation} we give an outline of some basic properties
of matrix groups; then in Section~\ref{sec:lfunctions} we give
some examples
of families of $L$-functions. In Section~\ref{sec:results}
we present the theorems for ratios of characteristic polynomials,
quoted from~\cite{CFZ} and~\cite{HPZ}.
In Section~\ref{sec:conjectures}
we describe our method of conjecturing precise formulas for averages of ratios
(an elaboration of the recipe in~\cite{CFKRS2}), and we give several
examples.
Refinements of the conjectures are presented in
Section~\ref{sec:refinements}. Finally, we mention some applications in
Section~\ref{sec:applications}.

\section{Random matrices and $L$-functions}\label{sec:notation}

We review the properties of the characteristic polynomials of
classical matrices
which we consider in this paper.

\subsection{Unitary matrices}  If $A=(a_{jk})$ is an
$N\times N$ matrix with complex entries, we let $A^*$
be its conjugate transpose, i.e. $A^*=(b_{jk})$ where
$b_{jk}=\overline{a_{kj}}$. $A$ is said to be unitary if $AA^*=I$.
We let $U(N)$ denote the group of all $N\times N$ unitary
matrices. This is a compact Lie group and has a Haar measure which
allows us to do analysis.

All of the eigenvalues of $A\in U(N)$ have absolute value~1; we write them as
 \begin{equation}e^{i\theta_1}, e^{i\theta_2}, \dots , e^{i\theta_N}\end{equation}
with
\begin{equation}0\le \theta_1 , \theta_2 , \dots , \theta_N < 2\pi.
\end{equation}
The eigenvalues of $A^*$ are $e^{-i\theta_1},\dots, e^{-i\theta_N}$.
The determinant, $\det A=\prod_{n=1}^N e^{i\theta_n}$,
is a complex number with absolute value~1.

For any sequence of $N$ points on the unit circle there are
matrices in $U(N)$ with those points as eigenvalues. The
collection of all matrices with the same set of eigenvalues
constitutes a conjugacy class in $U(N)$.  Thus, the set of
conjugacy classes can identified with the collection of sets
of $N$ points on the unit circle.


The characteristic polynomial of the matrix $A$ is denoted
$\Lambda_A(s)$ and is defined by
\begin{equation}\Lambda_A(s)=\det(I-sA^*)=\prod_{n=1}^N(1-se^{-i\theta_n}).
\end{equation}
The roots of $\Lambda_A(s)$ are the eigenvalues of $A$ and are on the unit circle.
The characteristic polynomial  satisfies the functional equation
\begin{align}
\nonumber
\Lambda_A(s)=\mathstrut&(-s)^N\prod_{n=1}^N e^{-i\theta_n}\prod_{n=1}^N
(1-e^{i\theta_n}/s)\cr
=\mathstrut&(-1)^N \det A^* ~s^N~\Lambda_{A^*}(1/s).
\end{align}
We regard $\Lambda_A(s)$ as an analogue of the Riemann zeta-function,
where the eigenangle $\theta$ plays the role of the parameter~$t$
on the critical line.

\subsection{Symplectic matrices}

The unitary matrix $A$ is said to be  {\it  symplectic} if
$AZA^t=Z$ where
\begin{equation*}
Z=\bigg(\begin{array}{cc} 0& I_N\\-I_N&0
\end{array}\bigg) ,
\end{equation*}
where $A^t$ denotes the transpose of~$A$.
A symplectic matrix necessarily has determinant~1.
The
symplectic group USp(2$N$) is the  group of $2N\times 2N$
unitary
symplectic matrices. The eigenvalues of a symplectic matrix
occur in complex conjugate pairs and we write them as
\begin{equation}
e^{\pm i\theta_1},\dots,e^{\pm i\theta_N}
\end{equation}

with
\begin{equation}
0\le\theta_1, \theta_2, \dots  ,\theta_N\le \pi.
\end{equation}

The functional equation of the characteristic polynomial is
\begin{equation}
\Lambda_A(s) = s^{2 N} ~\Lambda_{A^t}(1/s).
\end{equation}


\subsection{Orthogonal matrices}
A unitary matrix $A$ is said to be {\it orthogonal} if $AA^t=I$.
Orthogonality for a
unitary matrix implies that $A^t=A^*$ or $\overline{A}=A$. In
other words, an orthogonal matrix is a real unitary matrix.
We let $SO(N)$
denote the subgroup of $U(N)$ consisting of $N \times N$
orthogonal matrices with determinant~1.

The functional equation of the characteristic polynomial is
\begin{equation}
\Lambda_A(s) =(-1)^N ~ s^N ~  \Lambda_{A^t}(1/s).
\end{equation}
Thus, if $N$ is even  then the sign in the functional equation is 1
and if $N$ is odd then the sign is $-1$.
We distinguish these two cases as ``even'' orthogonal and ``odd'' orthogonal.

For any complex eigenvalue of an orthogonal matrix, its complex
conjugate is also an eigenvalue.  The eigenvalues of $A\in SO(2N)$
can be written as
\begin{equation}
e^{\pm i\theta_1},\dots,e^{\pm i\theta_N}
\end{equation}
and the eigenvalues of $A\in SO(2N+1)$ can be written as
\begin{equation}
1,e^{\pm i\theta_1},\dots,e^{\pm i\theta_{N }},
\end{equation}
where in both cases
\begin{equation}
0\le\theta_1, \theta_2,\dots  ,\theta_N\le \pi.
\end{equation}




\section{L-functions}\label{sec:lfunctions}

We briefly describe some of the families of $L$-functions for which we can
formulate a ratios conjecture.  See~\cite{F} for an introduction
and~\cite{CFKRS2} for all the details.  Much of the present paper is an
extension of~\cite{CFKRS2} and we assume some familiarity with that paper.

For the purpose of making conjectures
for ratios, the main feature of a family is that it is partially
ordered by a parameter (usually called the ``conductor''),
and there is a ``harmonic detector'' which describes the behavior
of the coefficients when averaged over the family.

\subsection{The Riemann zeta-function}
The Riemann zeta-function is defined by
\begin{equation}
\zeta(s)=\sum_{n=1}^\infty \frac{1}{n^s}
\end{equation}
for $s=\sigma+it$ with $\sigma >1.$ It has a
meromorphic continuation to the whole complex plane with its only
singularity a simple pole at $s=1$ with residue~1.  It satisfies a functional equation
which in its symmetric form reads
\begin{equation}
\pi^{-\frac s 2}\Gamma\bigg(\frac s 2 \bigg) \zeta(s) =
\pi^{ \frac {s-1} 2}\Gamma\bigg(\frac {1-s} 2 \bigg) \zeta(1-s)
\end{equation}
and in its asymmetric form
\begin{equation}
\zeta(s)=\chi(s)\zeta(1-s)
\end{equation}
where
\begin{equation}\label{eqn:fechi}
\chi(1-s)=\chi(s)^{-1}=2(2\pi)^{-s}\Gamma(s)\cos \frac{\pi s}{2}.
\end{equation}
The product formula discovered by Euler is
\begin{equation}
\zeta(s)=\prod_p \bigg( 1-\frac{1}{p^s}\bigg)^{-1}
\end{equation}
for $\sigma>1$ where the product is over the prime numbers $p$.

The family $\{\zeta(1/2+it)| t>0\}$ parametrized by real numbers $t$
can be modeled by characteristic polynomials of unitary matrices.
We will use a modification of the recipe in~\cite{CFKRS2} to conjecture
mean values for ratios of products of $\zeta$-functions.
The key ingredient in the recipe is the orthogonality
relation (or ``harmonic detector'')
\begin{equation}
\lim_{T\to \infty}\frac 1 T
\int_0^T \bigg(\frac m n\bigg)^{it}~dt = \bigg\{{ {1\mbox{ if $m=n$\phantom{.}}} \atop
{0\mbox{ if $m\ne n$.}}}
\end{equation}

\subsection{Dirichlet $L$-functions with real characters}
We let
\begin{equation}
L(s,\chi_d)=\sum_{n=1}^\infty \frac{\chi_d(n)}{n^s}
\end{equation}
for $\Re s>1$ where $\chi_d(n)$ is a primitive, real Dirichlet character.
The complete set of these characters is   described below.
Each of these (with $|d|>1$) is an entire function of $s$  and, if $d>0$,
 satisfies the functional equation
\begin{equation}
\bigg(\frac{\pi}{|d|}\bigg)^{-\frac s 2}\Gamma\bigg(\frac s 2 \bigg) L(s,\chi_d) =
\bigg(\frac{\pi}{|d|}\bigg)^{ \frac {s-1} 2}\Gamma\bigg(\frac {1-s} 2 \bigg) L(1-s,\chi_d)
\end{equation}
whereas if $d<0$, satisfies the functional equation
\begin{equation}
\bigg(\frac{\pi}{|d|}\bigg)^{-\frac s 2}\Gamma\bigg(\frac {s+1} 2 \bigg) L(s,\chi_d) =
\bigg(\frac{\pi}{|d|}\bigg)^{ \frac {s-1} 2}\Gamma\bigg(\frac {2-s} 2 \bigg) L(1-s,\chi_d).
\end{equation}

We now describe the characters $\chi_d$. These are
not defined for all $d$ but only for $d$ which are known as
fundamental discriminants. The values taken on by  $\chi_d(n)$
are $0$, $-1$ and $+1$.
We begin with
$\chi_{-4}(n)$ which is defined to be 1 if $n\equiv 1 \bmod 4$,   is
defined to be $-1$ if
$n\equiv 3 \bmod 4$, and is 0 if $n$ is even.
Next, we have
$\chi_{-8}(n)$ which is defined to be 0 if $n$ is even, is
$+1$ if $n\equiv 1$ or  $3 \bmod 8$ and is $-1$ if  $n\equiv 5$ or   $7\bmod 8$.
We also have
$\chi_8(n)$  which is defined to be 0 if $n$ is even, is
+1 if $n\equiv 1$ or  $7\bmod 8$ and is $-1$ if  $n\equiv 3$ or  $5\bmod 8$.
This takes care of all of the $d$ which are plus or minus a power of 2.
Now $d$ can also be equal to a prime $p \equiv 1 \bmod 4$.  In this
case $\chi_p(n)$ is 0 if $n$ is divisible by $p$, is
+1 if $n \equiv a^2 \bmod p$ for some $a$ not divisible by $p$ and
is $-1$ otherwise. If $p\equiv 3 \bmod 4$ then there is a character
$\chi_{-p}(n)$ which is defined exactly as
$\chi_p(n)$ for a $p\equiv 1 \bmod 4$. Finally, we
can take any pointwise product of distinct $\chi_d(n)$
to form
\begin{equation}
\chi_{d_1\cdot d_2 \cdot \dots \cdot d_u}(n)=
\chi_{d_1}(n)\chi_{d_2}(n)\dots \chi_{d_u}(n).
\end{equation}
If we take the empty product then we have the character
$\chi_1(n)=1$ for all $n$, so that $L(s,\chi_1)=\zeta(s)$.
This completes the description
of all of the primitive real characters.

Note that each $\chi_d(n)$ is defined for all integers $n$, positive and negative.
Also, $\chi_d(n)$ is periodic in $n$ with smallest period equal to $|d|$.
These functions are completely multiplicative, which means that
\begin{equation}
\chi_d(mn)=\chi_d(m)\chi_d(n).
\end{equation}
This multiplicativity implies
that $L(s,\chi_d)$ has an Euler product formula
\begin{equation}
L(s,\chi_d)=\prod_p\bigg(1-\frac{\chi_d(p)}{p^s}\bigg)^{-1}
\end{equation}
valid for $\Re s>1$.
If $d<0$ then $\chi_d$ is odd (i.e. $\chi_d(-n)=-\chi_d(n)$), whereas if
$d>0$, then $\chi_d(n)$ is even. We saw above that
the functional equations are slightly different in the even and odd
cases.

The collection of fundamental discriminants can be described as the set of
$d$ which   either are squarefree and congruent to 1 modulo 4 or are
4 times a squarefree number which is congruent to 2 or 3 modulo 4. The sequence
of $d$ is
\begin{equation}
\dots, -24,-23,-20,-19,-15,-11,-8,-7,-4,-3,1,5,8,12,13,17,21,24,28,29,33,37,\dots.
\end{equation}

The families $\{L(1/2,\chi_d)|d<0\}$ and  $\{L(1/2,\chi_d)|d>0\}$,
parametrized respectively by positive and negative fundamental discriminants,
can each be modeled by characteristic polynomials of symplectic matrices.
The harmonic detector for the positive discriminant family is
\begin{equation}
\delta(n):=\lim_{X \to \infty} \frac{1}{X^*} \sum_{0<d \le X}\chi_d(n)=
\bigg\{
{
{\prod_{p\mid n}(1+1/p)^{-1}    \mbox{ if $n$ is a perfect square}     }
\atop
{0     \mbox{           if $n$ is not a perfect square}     }
}
\end{equation}
where $X^*$ is the number of fundamental discriminants $0<d\le X$.
A similar formula holds for a sum over the odd characters $0<-d<X$.

\subsection{Quadratic twists of a modular $L$-function}
We now give an example of an orthogonal family of $L$-functions.

Let $f(z)$ be a holomorphic {\it newform} with integer
coefficients.  For example, if
\begin{equation}
   F(q)   =q\prod_{n=1}^\infty\left( 1-q^n\right) ^2\left(
1-q^{11n}\right) ^2  =\sum_{n=1}^\infty a_n q^n,
\end{equation}
where $a_1=1, a_2=-2, a_3=-1, a_4=2, \ldots$, then
\begin{equation*}
  f(z)=F\left( {\rm e}^{2\pi i z}\right) =
\sum_{n=1}^\infty a_n {\rm e}^{2\pi i n z}
\end{equation*}
is a newform (i.e.~it is an eigenfunction of the appropriate Hecke
operators). Specifically, $f(z)$ is a cusp form of weight 2 for
\begin{equation*}
 \Gamma_0(11)= \left\{\begin{pmatrix}
  a & b \\
  c & d
\end{pmatrix} \in SL(2,\mathbb{Z}) \ :\  11|c \right\};
\end{equation*}
that is
\begin{equation*}
\label{eq:4} f\left( \frac{az+b}{cz+d}\right)
=(cz+d)^2f(z)\;
\ \ \text{\ for all \ } \
 \begin{pmatrix}
  a & b \\
  c & d
\end{pmatrix} \in \Gamma_0(11).
\end{equation*}
The exponent of the factor multiplying $f(z)$  is the weight,
which is 2 in this case.

This newform is associated with
an elliptic curve.  Denote by $E_{11}$ the elliptic
curve
\begin{equation*}
y^2=4x^3-4x^2-40x-79,
\end{equation*}
and let
\begin{equation*}
 N_p=\#\{(x,y) :y^2\equiv 4x^3-4x^2-40x-79 \bmod p\},
\end{equation*}
then the coefficients    satisfy
\begin{equation*}
 a_p=p-N_p.
\end{equation*}

Deligne's Theorem  (proved earlier by Hasse for the special case of
elliptic curves) states that for a
weight-$k$ newform we have $|a_p|<2p^{\frac{k-1}{2}}$, and so in our
example $|a_p|<2\sqrt{p}$. We write
\begin{equation}
\lambda_E(n)=\lambda(n)=\frac{a(n)}{\sqrt{n}}.
\end{equation}
Therefore, one may associate with a newform $f$ a
Dirichlet series, called the $L$-function of the modular form,
\begin{equation}
L_{E_{11}}(s)   =\sum_{n=1}^\infty\frac{\lambda(n) }{n^{s }}
\end{equation}
which converges absolutely for $\Re s>1$.
The coefficients $\lambda(n)$ also satisfy the Hecke relations
\begin{equation}
\lambda(m)\lambda(n)=\sum_{d\mid m\atop d\mid n} \lambda\bigg(\frac{mn}{d^2}\bigg)
\end{equation}
which implies that the
$L$-function has an Euler product
\begin{equation}
L_{E_{11}}(s) = \left(1-\frac{1}{11^{s+1/2}}\right)^{-1}\prod_{p\neq11}
\left(1-\frac{\lambda(p)}{p^{s }}+\frac{1}{p^{2s}}\right)^{-1}.
\end{equation}
The $L$-function
associated with an elliptic curve is entire and satisfies the
functional equation
\begin{equation*}
 \left(\frac{2\pi}{\sqrt{M}}\right)^{-s}{\rm
\Gamma}(s)L_E(s)=w(E)\left(\frac{2\pi}{\sqrt{M}}\right)^{s-1}{\rm
\Gamma}(1-s)L_E(1-s),
\end{equation*}
where $M$ is the conductor of the elliptic curve $E$ and $w(E)=\pm
1$ is called the sign of the functional equation.  For $E_{11}$,
we have $M=11$ and $w(E)=1$.

The family we want to describe   is the collection of $L$-functions
associated with the quadratic twists of a fixed  $L$-function.  Let
$L_E$ be the $L$-function associated with an elliptic curve $E$
and let $\chi_d(n)$ be a real primitive Dirichlet character,
as described in the previous section.
Then the twisted
$L$-function
\begin{equation*}
L_{E}(s,\chi_d)=\sum_{n=1}^{\infty}\frac{\lambda(n)\chi_d(n)}{n^{s }}
\end{equation*}
is the $L$-function of another elliptic curve $E_d$: the quadratic
twist of $E$ by $d$. It can be shown that the number $N_{p,d}$ of solutions of
$E_d$ modulo $p$  satisfies
\begin{equation*}
 p-N_{p,d}=\chi_d(p)a_p.
\end{equation*}
Moreover,
$L_{E_d}(s)=L_{E}(s,\chi_d)$ satisfies the functional equation
\begin{equation*}
  \left(\frac{2\pi}{\sqrt{M}|d|}\right)^{-s}{\rm
\Gamma}(s)L_{E_d}(s)=\chi_d(-M)w(E)\left(\frac{2\pi}{\sqrt{M}|d|}\right)^{s-1}{\rm
\Gamma}(1-s)L_{E_d}(1-s).
\end{equation*}
For example, the quadratic twist of $E_{11}$ by $d$ is the
elliptic curve
\begin{equation*}
 dy^2=4x^3-4x^2-40x-79.
\end{equation*}
The corresponding twisted $L$-function is
\begin{align}
L_{E_{11}}(s,\chi_d)) & =\sum_{n=1}^\infty\left(\frac{d}{n}\right)\frac{\lambda(n)}
{n^{s }} \cr
&
=\left(1-\left(\frac{d}{11}\right)\frac{1}{11^{s+1/2}}\right)^{-1}\prod_{p\not{\,|} 11d}
\left(1-\left(\frac{d}{n}\right)\frac{\lambda(p)}{p^{s }}+\frac{1}{p^{2s}}\right)^{-1}.
\end{align}
This satisfies the functional equation
\begin{equation*}
 \left(\frac{2\pi}{\sqrt{11}|d|}\right)^{-s}{\rm
\Gamma}(s)L_{E_{11,d}}(s)=\chi_d(-11)\left(\frac{2\pi}{\sqrt{11}|d|}\right)^{s-1}{\rm
\Gamma}(1-s)L_{E_{11,d}}(1-s).
\end{equation*}

Note that when $d\equiv 2,6,7,8,$ or 10  $\bmod\ 11$, then the sign
in the functional equation is $+1$.  The $L_{E_{11}}(s,\chi_d)$ for these $d$
form an even orthogonal family. If
$d\equiv 1,3,4,5$ or 9 $\bmod\ 11$, then the sign
in the functional equation is $-1$ and the
$L_{E_{11}}(s,\chi_d)$ for these $d$
form an odd orthogonal family.

\section{Autocorrelation of ratios of characteristic polynomials}\label{sec:results}

We quote formulas from \cite{CFZ} and~\cite{HPZ} for the ratios of
characteristic polynomials averaged over the unitary, symplectic,
and orthogonal matrix groups.
Variants of these formulas have also been given by
Basor and Forrester~\cite{BF}, Day~\cite{D},
Baik, Deift, and Strahov~\cite{BDS}, Fyodorov and Strahov~\cite{FS},
and   others.
New proofs for these formulas have also recently been given by
Conrey, Forrester, and Snaith~\cite{CFS}.

Note that in the case of an equal number of characteristic polynomials
in the numerator and denominator,
the results we quote from \cite{CFZ} and~\cite{HPZ} are valid for
\emph{all}~$N$, while the other methods are only valid for
sufficiently large~$N$.

We let
\begin{equation}
z(x)=\frac{1}{1-e^{-x}}=\frac{1}{x}+O(1).
\end{equation}
It will be seen that the function $z(x)$ plays the role for random matrix
theory that
$\zeta(1+x)$ plays in the theory of moments of $L$-functions.

Also let $\Xi_{K,L}$ denote the subset of permutations $\sigma\in S_{K+L}$ of
$\{1,2,\dots,K+L\}$ for which
\begin{equation}
\sigma(1)<\sigma(2)<\dots < \sigma(K)
\end{equation}
and
\begin{equation}
\sigma(K+1)<\sigma(K+2)<\dots < \sigma(K+L).
\end{equation}
The cardinality of $\Xi_{K,L}$ is $\big({{K+L}\atop K}\big)=
\frac{(K+L)!}{K!L!}.$
Finally, let $\epsilon=(\epsilon_1,\dots,\epsilon_K)$ be a vector
with each component
$\epsilon_k=\pm 1$ and denote $\mbox{sgn}(\epsilon)=\prod_{k=1}^K \epsilon_k.$

We let $dA$ denote Haar measure on whichever group we are integrating over.

\begin{theorem}\label{thm:U}  If $N\ge \max\{Q-K, R-L\}$ and $\Re(\gamma_q),\Re(\delta_r)>0$
then
\begin{align}
\int_{U(N)}&\frac{\prod_{k=1}^K\Lambda_A(e^{-\alpha_j})\prod_{\ell=K+1}^{K+L}
\Lambda_{A^*}(e^{\alpha_\ell})} {\prod_{q=1}^Q\Lambda_A(e^{-\gamma_q})
\prod_{r= 1}^{R}
\Lambda_{A^*}(e^{-\delta_r})}dA\cr
&\qquad = \sum_{\sigma\in \Xi_{K,L}}
e^{N\sum_{k=1}^K (\alpha_{\sigma(k)}-\alpha_k)}\frac{
\prod_{k=1}^K \prod_{\ell=K+1}^{K+L}z(\alpha_{\sigma(k)}-
\alpha_{\sigma(\ell)})\prod_{q=1}^Q \prod_{r=1}^R
z(\gamma_q+\delta_r)}
{\prod_{r+1}^ R \prod_{ k=1}^K z( \alpha_{\sigma(k)}+\delta_r)
\prod_{q=1}^Q \prod_{ \ell=K+1}^{K+L} z(\gamma_q-\alpha_{\sigma(\ell)}) } .
\end{align}
\end{theorem}
If we let
\begin{equation}y_U(\alpha;\beta;\gamma;\delta):=
\frac{\prod_{k=1}^K\prod_{\ell=1}^L z( \alpha_k+\beta_\ell)\prod_{q=1}^Q
\prod_{r=1}^Rz( \gamma_q+\delta_r)}
{\prod_{k=1}^K\prod_{r=1}^R z( \alpha_k+\delta_r)
\prod_{\ell=1}^L\prod_{q=1}^Q z( \beta_\ell+\gamma_q)},\end{equation}
then the above can be expressed as
\begin{align}
\int_{U(N)}&\frac{\prod_{k=1}^K\Lambda_A(e^{-\alpha_k})\prod_{\ell=K+1}^{K+L}
\Lambda_{A^*}(e^{\alpha_\ell})} {\prod_{q=1}^Q\Lambda_A(e^{-\gamma_q})
\prod_{r= 1}^{R}
\Lambda_{A^*}(e^{-\delta_r})}dA
\cr
&\qquad =\sum_{\sigma \in \Xi_{K,L}}
e^{N\sum_{k=1}^K (\alpha_{\sigma(k)}-\alpha_k)}
y_U(\alpha_{\sigma(1)},\dots,  \alpha_{\sigma(K)};
-\alpha_{\sigma(K+1)}\dots  -\alpha_{\sigma(K+L)};\gamma;\delta) .
\end{align}

\begin{theorem}\label{thm:Sp}    If $2N \ge Q-K-1$ and $\Re(\gamma_q)>0$
then
\begin{align}
\int_{USp(2N)}&\frac{\prod_{k=1}^K\Lambda_A(e^{-\alpha_k})}
{\prod_{q=1}^Q \Lambda_A(e^{-\gamma_q})}~dA\cr
&\qquad =
\sum_{\epsilon \in \{-1,1\}^K} e^{N\sum_{k=1}^K(\epsilon_k\alpha_k-\alpha_k)}
\frac{\prod_{  j \le k\le K} z(\epsilon_j \alpha_j + \epsilon_k \alpha_k)
\prod_{  q < r\le Q} z( \gamma_q +   \gamma_r)}
{
\prod_{k=1}^K  \prod_{q=1}^Q  z(\epsilon_k \alpha_k +   \gamma_q)},
\end{align}
\end{theorem}

If we let
\begin{equation}y_S(\alpha ;\gamma ):=
\frac{\prod_{  j\le k \le K}z( \alpha_j+\alpha_k)\prod_{  q < r \le Q}
z( \gamma_q+\gamma_r)}
{\prod_{k=1}^K\prod_{q=1}^Q z( \alpha_k+\gamma_q)
}\end{equation}
and
\begin{equation}
h_S(\alpha;\gamma)= e^{N\sum_{k=1}^K\epsilon_k\alpha_k }y_S(\alpha;\gamma),
\end{equation}
then the above can be expressed as
\begin{equation}
 \int_{USp(2N)}\frac{\prod_{k=1}^K\Lambda_A(e^{-\alpha_k})}
{\prod_{q=1}^Q \Lambda_A(e^{-\gamma_q})}~dA
  =e^{-N\sum_{k=1}^K \alpha_k}
\sum_{\epsilon \in \{-1,1\}^K} h_S(\epsilon_1 \alpha_{1},\dots, \epsilon_K  \alpha_{K};
 \gamma ) .
\end{equation}

\begin{theorem}\label{thm:Oe}     If $2N \ge Q-K+1$ and $\Re(\gamma_q)>0$
then
\begin{align}
\int_{SO(2N)}&\frac{\prod_{k=1}^K\Lambda_A(e^{-\alpha_k})}
{\prod_{q=1}^Q \Lambda_A(e^{-\gamma_q})}~dA\cr
&\qquad =
\sum_{\epsilon \in \{-1,1\}^K} e^{N\sum_{k=1}^K(\epsilon_k\alpha_k-\alpha_k)}
\frac{\prod_{  j < k\le K} z(\epsilon_j \alpha_j + \epsilon_k \alpha_k)
\prod_{  q \le r\le Q} z( \gamma_q +   \gamma_r)}
 {\prod_{k=1}^K \prod_{r=1}^R    z(\epsilon_k \alpha_k +   \gamma_r)
}.
\end{align}
\end{theorem}
If we let
\begin{equation}y_O(\alpha ;\gamma ):=
\frac{\prod_{  j< k \le K}z( \alpha_j+\alpha_k)\prod_{  q \le r \le Q}
 z( \gamma_q+\gamma_r) }
{\prod_{k=1}^K\prod_{q=1}^Q z( \alpha_k+\gamma_q)
 }\end{equation}
and
\begin{equation}
h_O(\alpha;\gamma)= e^{N\sum_{k=1}^K\epsilon_k\alpha_k }y_O(\alpha;\gamma),
\end{equation}
then the above can be expressed as
\begin{equation}
 \int_{SO(2N)}\frac{\prod_{k=1}^K\Lambda_A(e^{-\alpha_k})}
{\prod_{q=1}^Q \Lambda_A(e^{-\gamma_q})}~dA
  =
e^{-N\sum_{k=1}^K \alpha_k}
\sum_{\epsilon \in \{-1,1\}^K} h_O(\epsilon_1 \alpha_{1},\dots, \epsilon_K  \alpha_{K};
 \gamma )  .
\end{equation}

\begin{theorem}\label{thm:Oo} If $2N\ge Q-K$ and $\Re(\gamma_q)>0$
\begin{align}
\int_{SO(2N+1)}&\frac{\prod_{k=1}^K\Lambda_A(e^{-\alpha_k})}
{\prod_{q=1}^Q \Lambda_A(e^{-\gamma_q})}~dA\cr
&\qquad =
\sum_{\epsilon \in \{-1,1\}^K}\mbox{sgn}(\epsilon) e^{N\sum_{k=1}^K
(\epsilon_k\alpha_k-\alpha_k)}
\frac{\prod_{  j < k\le K} z(\epsilon_j \alpha_j + \epsilon_k \alpha_k)
\prod_{ q \le r\le Q} z( \gamma_q +   \gamma_r)}
 {\prod_{k=1}^K\prod_{q=1}^Q  z(\epsilon_k \alpha_k+   \gamma_q)
}\cr
& \qquad=
e^{-N\sum_{k=1}^K \alpha_k}
\sum_{\epsilon \in \{-1,1\}^K}
\mbox{ sgn}(\epsilon) h_O(\epsilon_1 \alpha_{1},\dots, \epsilon_K  \alpha_{K};
 \gamma ).
\end{align}
\end{theorem}

In the next section we give conjectures for the averages of ratios
of $L$-functions, which will have a very similar form to the theorems
given above.

\section{Conjectures about autocorrelations of ratios of $L$-functions}\label{sec:conjectures}

We  make conjectures about averages of $L$-functions which are analogous
to the theorems of the previous section. Roughly speaking, the number
$N$ of independent eigenvalues
of the
matrix   is replaced
by the analytic conductor of the family,
the function $z(x)$ in the above  theorems
is replaced by $\zeta(1+x)$, and in addition, an arithmetic factor
$A$   must be introduced.
This arithmetic factor, which depends on the
particular family under consideration,
is expressible as an infinite
product over primes and can be computed
on a case by case basis.

In Section~\ref{sec:recipe} we give our recipe for conjecturing averages
of ratios, then we illustrate the computation for
some standard examples of families of $L$-functions.
For a more detailed discussion of related conjectures, see~\cite{CFKRS2}.

\subsection{The recipe}\label{sec:recipe}

The following is an extension of
the approximate functional equation recipe of~\cite{CFKRS2}.
Familiarity with that paper will be helpful here.

Suppose $\mathcal L$ is an $L$-function and ${\mathcal F} = \{f\}$ is a
family of characters with conductor~$c(f)$,
as described in Section~3 of~\cite{CFKRS2}.
Thus, $\mathcal{L}(s,f)$ has an approximate functional
equation of the form
\begin{equation}\label{eqn:generalappfe}
\mathcal{L}(s,f)= \sum \frac{a_n(f)}{n^s}+
\varepsilon_f \mathcal{X}_f(s) \sum
\frac{\overline{a_n(f)}}{n^{1-s}} + remainder.
\end{equation}
Also, we can write
\begin{equation}\label{eqn:reciprocal}
\frac{1}{\mathcal{L}(s,f)} =
\sum_{n=1}^\infty
\frac{\mu_{\mathcal{L},f}(n)}{n^s},
\end{equation}
the series converging absolutely for $\Re(s)>1$ and conditionally,
assuming a suitable Riemann Hypothesis, for $\Re(s)>\frac12$.

We wish to conjecture a precise asymptotic formula for the average
\begin{equation}\label{eqn:generalmoment}
\sum_{f\in \mathcal F}
\frac{\mathcal{L}(\tfrac12 +\alpha_1,f)\dots
\mathcal{L}(\tfrac12 +\alpha_K,f)
\mathcal{L}(\tfrac12 +\alpha_{K+1},\overline{f})\dots
\mathcal{L}(\tfrac12 +\alpha_{K+L},\overline{f})}
{
\mathcal{L}(\tfrac12 +\gamma_1,f)\dots
\mathcal{L}(\tfrac12 +\gamma_Q,f) \mathcal{L}(\tfrac12 +\delta_{1},\overline{f})\dots
\mathcal{L}(\tfrac12 +\delta_{R},\overline{f})}
g(c(f))
\end{equation}
where $g$ is a suitable test function.
Note that the sum is an integral in the case of moments in $t$-aspect.

The recipe:

\begin{enumerate}
\item{} Start with
\begin{align}
\mathcal{L}_f(s; \boldsymbol\alpha_K;\boldsymbol\alpha_{L};
        \boldsymbol\gamma_Q; \boldsymbol\delta_{R}) =&
\frac{\mathcal{L}(s +\alpha_1,f)\dots
\mathcal{L}(s +\alpha_K,f)
\mathcal{L}(s +\alpha_{K+1},\overline{f})\dots
\mathcal{L}(s +\alpha_{K+L},\overline{f})}
{
\mathcal{L}(s +\gamma_1,f)\dots
\mathcal{L}(s +\gamma_Q,f) \mathcal{L}(s +\delta_{1},\overline{f})\dots
\mathcal{L}(s +\delta_{R},\overline{f})}
\end{align}

\item{} Replace each $L$-function in the numerator with the two terms from its
approximate functional equation \eqref{eqn:generalappfe}, ignoring the remainder term.
Replace each $L$-function in the denominator by its series~\eqref{eqn:reciprocal}.
Multiply out the resulting
expression to obtain $2^{K+L}$ terms.  Write those terms as
\begin{equation}
(\text{product of $\varepsilon_f$ factors})
(\text{product of $\mathcal{X}_f$ factors})
\sum_{n_1,\dots,n_{K+L+Q+R}} (\text{summand}) .
\end{equation}

\item{} Replace each product of $\varepsilon_f$-factors
by its expected value when averaged over the family.

\item{} Replace each summand by its expected value when averaged over the family.

\item{} Complete the resulting sums (i.e., extend the ranges of the
summation indices out to infinity), and call the total
$M_f(s,\boldsymbol\alpha_K;\boldsymbol\alpha_{L};
        \boldsymbol\gamma_Q; \boldsymbol\delta_{R})$.

\item{} The conjecture is
\begin{align}
\sum_{f\in \mathcal F} {\mathcal L}_f(\tfrac12,
 \boldsymbol\alpha_K;\boldsymbol\alpha_{L};
        \boldsymbol\gamma_Q; \boldsymbol\delta_{R}
        ) g(c(f)) 
=&
\sum_{f\in \mathcal F} M_f(\tfrac12, \boldsymbol\alpha_K;\boldsymbol\alpha_{L};
        \boldsymbol\gamma_Q; \boldsymbol\delta_{R})
(1+O(e^{(-\frac12 + \varepsilon)c(f)})) g(c(f)),
\end{align}
for all $\varepsilon>0$, where $g$ is a suitable weight function.
\end{enumerate}

In other words, ${\mathcal L}_f(\tfrac12,\cdot)$ and $M_f(\tfrac12,\cdot)$ have the same value
distribution if averaged over a sufficiently large portion of the family.
Note that the dependence of $M_f$ on $f$ only occurs in the product
of $\mathcal{X}_f$ factors.

The above conjecture has a square-root error term.  Presumably this
is best possible.  Since very little is known about mean-values of ratios,
we are not able to give any objective evidence for such
a small error term.  Also, we have not specified the allowable
range for the shifts $\alpha$, $\gamma$, and $\delta$.  Conrey and
Snaith~\cite{CS} suggest that in the case of the zeta-function
one should allow shifts with imaginary part~$\ll T$.

\subsection{Moments of ratios of $\zeta(s)$}
Let $s=1/2+it$ and consider
\begin{equation} \frac{1}{T}\int_0^T
\frac{\prod_{k=1}^{K}\zeta(s+\alpha_k)\prod_{\ell=K+1}^{K+L}
\zeta(1-s-\alpha_\ell)}{\prod_{q=1}^Q\zeta(s+\gamma_q)
\prod_{r=1}^R\zeta(1-s+\delta_r)} ~dt\end{equation}
where $\Re \alpha_k,\gamma_q,\delta_r > 0$.
We assume the Riemann Hypothesis so there are no poles
on the path of integration.
We follow the recipe given above.

First replace each $\zeta$-function in the numerator by its approximate
 functional equation
\begin{equation}\zeta(s)\sim \sum_{n\le \tau} \frac{1}{n^s}+\chi(s)
\sum_{n\le \tau}\frac{1}{n^{1-s}}\end{equation}
(where $\tau=\sqrt{t/(2\pi)}$).
Here $\chi$ is the factor in the functional equation, given in
\eqref{eqn:fechi}.
 Second, multiply out to get $2^{K+L}$ terms.
Only retain the   terms in which the same number of $\chi(s)$
as $\chi(1-s)$ occur, because the others are highly oscillatory and have
expected value 0 when averaged over~$t$.
For the terms in the denominator,
expand
into series
\begin{equation}\frac{1}{\zeta(s)}=\prod_p\bigg(1-\frac{1}{p^s}\bigg)=
\sum_{n=1}^\infty \frac{\mu(n)}{n^s}.
\end{equation}
Here $\mu(n)$ is the M\"{o}bius function, which is multiplicative and
is equal to $-1$ when $n=p$ is prime and is 0 when $n=p^e$ where $e>1$.

For each of the retained terms, we keep those summands which
are independent of the parameter~$t$; equivalently, we
keep the ``diagonal.''
(The other summands are of the form $\Theta^t$ with $\Theta\neq 1$,
which has 0 expected value.)

In these calculations it is easiest to initially work with
the
expression obtained from the `first' term of each approximate
functional equation, and then modify that expression to obtain the
complete main term.
This will determine the arithmetic factor $A$ which appears in the
$L$-function averages (but not in the analogous random matrix averages).
For convenience, let $\beta_\ell=-\alpha_{K+\ell}$; moreover, assume that the real parts
of all the variables $\alpha_k,\beta_\ell,\gamma_q, \delta_r$ are positive
so that the series
\begin{equation}G_\zeta(\alpha;\beta; \gamma;\delta)=\sum_{\prod
m_k \prod h_q=\prod  n_\ell \prod j_r}
\frac{\prod \mu(h_q)\prod \mu(j_r)}
{\prod   m_k^{1/2+\alpha_k}\prod n_\ell^{1/2+\beta_\ell}\prod h_q^
{1/2+\gamma_q}\prod j_r^{1/2+\delta_r}}
\end{equation}
is absolutely convergent.
We  express $G_\zeta$ as an Euler product
\begin{equation}G_\zeta(\alpha;\beta;\gamma;\delta)=\prod_p \sum_{\sum a_k+\sum c_q=\sum b_\ell+
\sum d_r}
 \frac{\prod  \mu(p^{c_q})\prod \mu(p^{d_r})}{\prod p^{\sum_k(1/2+\alpha_k)
 a_k+\sum_\ell (1/2+\beta_\ell)b_\ell+
 \sum_q(1/2+\gamma_q)c_q+\sum_r(1/2+\delta_r)d_r}} .
\end{equation}
The terms here with $\sum_{k=1}^K a_k+\sum_{q=1}^Q c_q=1=\sum_{\ell=1}^n b_\ell
+\sum_{r=1}^R d_r$
 contribute the zeros and poles.   These terms give
\begin{equation}Y_U(\alpha;\beta;\gamma;\delta):=
\frac{\prod_{k=1}^K\prod_{\ell=1}^L \zeta(1+\alpha_k+\beta_\ell)\prod_{q=1}^Q
\prod_{r=1}^R\zeta(1+\gamma_q+\delta_r)}
{\prod_{k=1}^K\prod_{r=1}^R \zeta(1+\alpha_k+\delta_r)
\prod_{\ell=1}^L\prod_{q=1}^Q \zeta(1+\beta_\ell+\gamma_q)}.\end{equation}
 We factor $Y_U$ out of $G_\zeta$ and are left with
\begin{equation}G_\zeta(\alpha;\beta;\gamma;\delta)=Y_U(\alpha;\beta;\gamma;\delta)
A_\zeta(\alpha;\beta;\gamma;\delta)\end{equation}
where
$A_\zeta$ is an Euler product,  absolutely convergent
for all of the variables in small disks around 0, which is given by
\begin{align}
A_\zeta=&\prod_p \frac{\prod_{K=1 }^K\prod_{\ell=1}^L (1-1/p^{1+\alpha_k+\beta_\ell})
\prod_{q=1}^Q\prod_{r=1}^R(1-1/p^
{1+\gamma_q+\delta_r})}{\prod_{k=1}^K\prod_{r=1}^R(1-1/p^{1+\alpha_k+\delta_r})
\prod_{\ell=1}^L\prod_{q=1}^Q(1-1/p^{1+\beta_\ell+\gamma_q}) }\cr
& \qquad \times
\sum_{\sum  a_k+\sum c_q=\sum b_\ell+
\sum d_r}
\frac{\prod  \mu(p^c_q)\prod \mu(p^{d_r})}{  p^{\sum(1/2+\alpha_k)
a_k+\sum (1/2+\beta_\ell)b_\ell+\sum(1/2+\gamma_q)c_q+\sum(1/2+\delta_r)d_r}}
\end{align}

Thus, the recipe leads us to:
\begin{conjecture}\label{conj:zeta}  If $\Re(\gamma_q),\Re(\delta_r)>0$
and $\Re(\alpha_j)>-\frac{1}{2(K+L)}$, then
\begin{align}
\int_0^T &
\frac{\prod_{k=1}^{K}\zeta(s+\alpha_k)\prod_{\ell=K+1}^{K+L}
\zeta(1-s-\alpha_\ell)}{\prod_{q=1}^Q\zeta(s+\gamma_q)
\prod_{r=1}^R\zeta(1-s+\delta_r)} ~dt\cr
& = \int_0^T
\sum_{\sigma \in \Xi_{K,L}}\prod_{k=1}^K\frac{\chi(s+\alpha_k)}{\chi(s-\alpha_{\sigma(k)})}
Y_UA_\zeta(\alpha_{\sigma(1)},\dots,  \alpha_{\sigma(K)};
-\alpha_{\sigma(K+1)}\dots  -\alpha_{\sigma(K+L)};\gamma;\delta) ~dt\cr
& \qquad \qquad
+O(T^{1/2+\epsilon})
\end{align}
\end{conjecture}
Note that $\prod_{k=1}^{K+L}\chi(s+\alpha_k)=\prod_{k=1}^{K+L}\chi(s+\alpha_{\sigma(k)})$ so that
\begin{equation}
\prod_{k=1}^K\frac{\chi(s+\alpha_k)}{\chi(s-\alpha_{\sigma(k)})}=\prod_{k=1}^K\frac{\chi(s+\alpha_k)^{ 1/2}}{\chi(s-\alpha_{\sigma(k)})^{1/2}}
\prod_{\ell=1}^L\frac{\chi(s-\alpha_{\sigma(K+\ell)})^{1/2}}{\chi(s+\alpha_{K+\ell})^{1/2}}
\end{equation}
Thus, the factor on the left can be replaced by the factor on the right;
this leads to a slightly different formulation of the conjecture
which is convenient in Section~\ref{sec:combsum} where we replace the combinatorial sum by a
multiple integral.  In particular, letting
\begin{equation}
H_{\zeta,t}(w_1, \dots,w_{K+L};\gamma;\delta)=
\frac{\prod_{\ell=1}^L\chi(s-w_{K+\ell})^{1/2}}{\prod_{k=1}^K \chi(s+w_k)^{1/2}}
Y_UA_\zeta(w_1,\dots,w_K; -w_{K+1},\dots ,-w_{K+L};\gamma;\delta)
\end{equation}
the conjecture may be reformulated as
\begin{align}
\int_0^T &
\frac{\prod_{k=1}^{K}\zeta(s+\alpha_k)\prod_{\ell=K+1}^{K+L}
\zeta(1-s-\alpha_\ell)}{\prod_{q=1}^Q\zeta(s+\gamma_q)
\prod_{r=1}^R\zeta(1-s+\delta_r)} ~dt\cr
&= \int_0^T
\frac{\prod_{k=1}^K\chi(s+\alpha_k)^{1/2}}{\prod_{\ell=1}^L \chi(s-\alpha_{K+\ell})^{1/2}}
\sum_{\sigma \in \Xi_{K,L}}
 H_{\zeta,t}(\alpha_{\sigma(1)},\dots,  \alpha_{\sigma(K)};
 \alpha_{\sigma(K+1)}\dots   \alpha_{\sigma(K+L)};\gamma;\delta) ~dt\cr
& \qquad \qquad
+O(T^{1/2+\epsilon})
\end{align}

Conrey and Snaith~\cite{CS} have an alternative formulation for this conjecture with a subscript-free notation.

\subsection{Moments of $L(1/2,\chi_d)$}
The family $\mathcal D^+=\{L(s,\chi_d):d>0\}$ is a symplectic family. We can make a conjecture
analogous to Theorem~\ref{thm:Sp} for
\begin{equation}
\sum_{0<d\le X} \frac{\prod_{k=1}^K L(1/2+\alpha_k,\chi_d)}
{\prod_{q=1}^Q L(1/2+\gamma_q,\chi_d)} .
\end{equation}
As in the previous example, the main issue will be
identifying the appropriate arithmetic
factor $A_\mathcal D$.

Again we follow the recipe, which will parallel the example of the
Riemann $\zeta$-function in the previous section.
The $L$-functions in the numerator
are replaced by their approximate functional equations
while the ones in the denominator
are expanded into series
\begin{equation}
\frac{1}{L(s,\chi_d)}=\prod_p\bigg(1-\frac{\chi_d(p)}{p^s}\bigg)
=\sum_{n=1}^\infty \frac{\mu(n)\chi_d(n)}{n^s}
\end{equation}
with $\mu(n)$ as before.

Expanding the product of approximate functional equations,
we obtain $2^K$ terms.  All those terms are retained because
the sign of the functional equation is always~$+1$.
So now we replace the summands by their average.

We can  determine $G_\mathcal D$ (analogous to $G_\zeta$ in the
previous example) by consideration of
\begin{align}
G_\mathcal D(\alpha;\gamma):=&\lim_{X\to \infty} \frac{1}{X^*}
\sum_{m_k, h_q} \frac{\prod_q \mu(h_q)}
{\prod_k m_k^{1/2+\alpha_k} \prod_q h_q^{1/2+\gamma_q}}
\sum_{0<d\le X} \chi_d\big(\prod_k m_k \prod_q h_q\big)\cr
  =&
\sum_{m_k, h_q} \frac{\prod_q \mu(h_q)}
{\prod_k m_k^{1/2+\alpha_k} \prod_q h_q^{1/2+\gamma_q}}
  \delta\big(\prod_k m_k \prod_q h_q\big)
\end{align}
where $\delta(n)=\prod_{p\mid n}(1+1/p)^{-1}$ if $n$ is a square and is 0 otherwise.
We can express $G_\mathcal D$ as a convergent Euler product provided that
the real parts of the $\alpha_k$ and the $\gamma_q$ are positive.
Thus,
\begin{equation}
G_\mathcal D(\alpha;\gamma)=\prod_p \bigg(1+(1+\tfrac 1p)^{-1}
\sum_{0< \sum_k a_k +\sum_q c_q  \mbox{ is even}}
\frac{\prod_q \mu(p^{c_q})}
{  p ^{\sum_k a_k(1/2+\alpha_k) +\sum_q c_q(1/2+\gamma_q)}}\bigg) .
\end{equation}
The terms here with $\sum_{k=1}^K a_k+\sum_{q=1}^Q c_q=2 $
contribute the zeros and poles. Specifically, poles arise from
terms $a_j=a_k=1$ with $1\le j< k\le K$ and from terms
$a_k=2$ with $1\le k\le K$. Poles also arise from
terms with $c_q=c_r =1$ with $1\le q< r\le Q$.  Note that poles
do not arise from terms with $c_q=2$ since $\mu(p^2)=0$. Zeros arise
from terms with $a_k=1=c_q$ with $1\le k\le K$ and $1\le q\le Q$.    The  contribution,
expressed in terms of zeta-functions, of all of these zero and polar terms is
\begin{equation}Y_S(\alpha; \gamma ):=
\frac{\prod_{  j\le k\le K}  \zeta(1+\alpha_j+\alpha_k)\prod_{  q< r\le Q}
\zeta(1+\gamma_q+\gamma_r)}
{\prod_{k=1}^K\prod_{q=1}^Q \zeta(1+\alpha_k+\gamma_q)}.
\end{equation}
When we factor $Y_S$ out from $G_\mathcal D$ we are left
with an Euler product $A_\mathcal D(\alpha,\gamma)$
which is absolutely convergent for
all of the variables in small disks around 0. Specifically,
\begin{align}
A_\mathcal D(\alpha,\gamma)
 =&\prod_p \frac{\prod_{  j\le k\le K}  (1-1/p^{1+\alpha_j+\alpha_k})
\prod_{ q< r\le Q}(1-1/p^
{1+\gamma_q+\gamma_r})}{\prod_{k=1}^K\prod_{q=1}^Q(1-1/p^{1+\alpha_k+\gamma_q})
  }\cr
& \qquad \times
\bigg(1+(1+\tfrac 1p)^{-1}
\sum_{0< \sum_k a_k +\sum_q c_q  \mbox{ is even}}
\frac{\prod_q \mu(p^{c_q})}
{  p ^{\sum_k a_k(1/2+\alpha_k) +\sum_q c_q(1/2+\gamma_q)}}\bigg) .
\end{align}

The functional equation may be written as
\begin{equation}
L(s,\chi_d)=\bigg(\frac{|d|}{\pi}\bigg)^{ \frac 12-s}g_+(s) L(1-s,\chi_d)
\end{equation}
where
\begin{equation}
g_+(s)=\frac{\Gamma\big(\frac{1-s}{2}\big)}{\Gamma\big(\frac{ s}{2}\big)}.
\end{equation}
Note that $g_+(1/2)=1$.
The analytic conductor of $L(s,\chi_d)$ is $\frac{|d|}{ \pi}$ so that
the role of $2N$ in Theorem~\ref{thm:Sp} is played by $2N=\log \frac{|d|}{ \pi}$.
(There are some subtleties concerning the ``conductor'' here.
See the discussion of Conjecture~1.5.3 in~\cite{CFKRS2}).
We are led to
\begin{conjecture}\label{conj:realchi}
Suppose that the real parts of $\alpha_k$ and
$\gamma_q$ are positive. Then
\begin{align}
\sum_{0<d\le X} &\frac{\prod_{k=1}^K L(1/2+\alpha_k,\chi_d)}
{\prod_{q=1}^Q L(1/2+\gamma_q,\chi_d)}\cr
&\qquad  =
\sum_{0<d\le X} \sum_{\epsilon \in \{-1,1\}^K}
\bigg(\frac{|d|}{\pi}\bigg)^{\frac 12\sum_{k=1}^K(\epsilon_k\alpha_k-\alpha_k)}
\prod_{k=1}^K g_+\left(\frac 12+\frac{\alpha_k-\epsilon_k \alpha_k }{2}\right)
Y_SA_\mathcal D(\epsilon_1 \alpha_{1},\dots, \epsilon_K  \alpha_{K};
 \gamma )\cr
 &\qquad \qquad  +O(X^{1/2+\epsilon}) .
\end{align}
\end{conjecture}
If we let
\begin{equation}
H_{\mathcal D^+,d,\alpha,\gamma}(w)=
\bigg(\frac{|d|}{\pi}\bigg)^{\frac 12\sum_{k=1}^Kw_k}
\prod_{k=1}^K g_+\left(\frac 12+\frac{\alpha_k-w_k }{2}\right)
Y_SA_\mathcal D(w_1,\dots,w_k;
 \gamma )
\end{equation}
then
the conjecture may be formulated as
\begin{align}
\sum_{0<d\le X} &\frac{\prod_{k=1}^K L(1/2+\alpha_k,\chi_d)}
{\prod_{q=1}^Q L(1/2+\gamma_q,\chi_d)}\cr
&\qquad  =
\sum_{0<d\le X}
\bigg(\frac{|d|}{\pi}\bigg)^{-\frac 12\sum_{k=1}^K \alpha_k}
\sum_{\epsilon \in \{-1,1\}^K}
H_{\mathcal D^+,d,\alpha,\gamma}(\epsilon_1\alpha_1,\dots,\epsilon_K \alpha_K)  +O(X^{1/2+\epsilon}) .
\end{align}

A conjecture for the moments of ratios of $L(s,\chi_d)$ with $d<0$ can be
analogously formulated; the only change is that $g_+$ is replaced by
\begin{equation} g_-(s)=\frac{\Gamma\big(\frac{2-s}{2}\big)}{\Gamma\big(\frac{s+1}{2}\big)}.
\end{equation}

\subsection{Moments of $L_E(1/2,\chi_d)$}
Given an elliptic curve $E$,
the family $E^+(\mathcal D)=\{L_E(s,\chi_d):w(E_d)=+1 \}$ is an even orthogonal family
and   $E^-(\mathcal D)=\{L_E(s,\chi_d):w(E_d)=-1 \}$ is an odd orthogonal family.
We can  formulate   conjectures
analogous to Theorems~\ref{thm:Oe} and~\ref{thm:Oo} for
\begin{equation}
\sum_{|d|\le X\atop w(E_d)=+1} \frac{\prod_{k=1}^K L_E(1/2+\alpha_k,\chi_d)}
{\prod_{q=1}^Q L_E(1/2+\gamma_q,\chi_d)},
\end{equation}
and  for the same sum over $d$ with $w(E_d)=-1,$
once we have identified the appropriate arithmetic
factor $A_{E(\mathcal D)}$. Suppose that the $L$-function associated with $E$ has level $M$.   Let
\begin{equation}
\frac{1}{L_E(s)}=\prod_p\bigg(1-\frac{\lambda(p)}{p^s}+\frac{\chi_0(p)}{p^{2s}}\bigg)
=\sum_{n=1}^\infty \frac{\mu_E(n)}{n^s},
\end{equation}
where $\chi_0$ is the principal character $\bmod M$.
So $\mu_E(n)$ is a multiplicative function which is equal
to $-\lambda(p)$ for $n=p$, is equal to $\chi_0(p)$ if $n=p^2$, and is 0 if $n=p^e$ with $e>2$.
Let $X^*_+=|\{d:|d|\le X, w(E_d)=1\}|.$

As in the previous two examples, the calculation is fairly straightforward
up to the point of computing the arithmetic factor.  So
we consider
\begin{align}
G_{E(\mathcal D)}(\alpha;\gamma):=&\lim_{X\to \infty} \frac{1}{X_+^*}
\sum_{m_k, h_q} \frac{\prod_k \lambda(m_k) \prod_q \mu_E(h_q)}
{\prod_k m_k^{1/2+\alpha_k} \prod_q h_q^{1/2+\gamma_q}}
\sum_{|d|\le X\atop w(E_d)=1} \chi_d\big(\prod_k m_k \prod_q h_q\big)\cr
  =&
\sum_{m_k, h_q} \frac{\prod_k \lambda(m_k) \prod_q \mu_E(h_q)}
{\prod_k m_k^{1/2+\alpha_k} \prod_q h_q^{1/2+\gamma_q}}
  \delta\big(\prod_k m_k \prod_q h_q\big)
\end{align}
where, as before,
 $\delta(n)=\prod_{p\mid n}(1+1/p)^{-1}$ if $n$ is a square and is 0 otherwise.
We can express $G_{E(\mathcal D)}$ as a convergent Euler product provided that
the real parts of the $\alpha_k$ and the $\gamma_q$ are positive.
Thus,
\begin{equation}
G_{E(\mathcal D)}(\alpha;\gamma)=\prod_p \bigg(1+(1+\tfrac 1p)^{-1}
\sum_{0< \sum_k a_k +\sum_q c_q  \mbox{ is even}}
\frac{\prod_k \lambda(p^{a_k})\prod_q \mu_E(p^{c_q})}
{  p ^{\sum_k a_k(1/2+\alpha_k) +\sum_q c_q(1/2+\gamma_q)}}\bigg) .
\end{equation}
The terms with the smallest positive exponents, that is $\sum_{k=1}^K a_k+\sum_{q=1}^Q c_q=2 $,
contribute the zeros and poles.
Specifically, poles arise from
terms $a_j=a_k=1$ with $1\le j< k\le K$. Note that the  terms with
$a_k=2$ do not contribute poles; this is because the function
\begin{equation}
\sum_{n=1}^\infty \frac{\lambda(n^2)}{n^s}
\end{equation}
is analytic at $s=1$.
 Poles also arise from
terms with $c_q=c_r =1$ with $1\le q< r\le Q$.  Zeros arise
from terms with $a_k=1=c_q$ with $1\le k\le K$ and $1\le q\le Q$.
 Zeros also arise from terms with $c_q=2$ since $\mu_E(p^2)=1$.
  The  contribution,
expressed in terms of zeta-functions, of all of these zero and polar terms is
 \begin{equation}Y_O(\alpha; \gamma ):=
\frac{\prod_{  j< k\le K}  \zeta(1+\alpha_j+\alpha_k)\prod_{  q< r\le Q}
\zeta(1+\gamma_q+\gamma_r)\prod_{q=1}^Q \zeta(1+2\gamma_q)}
{\prod_{k=1}^K\prod_{q=1}^Q \zeta(1+\alpha_k+\gamma_q)}.
\end{equation}
When we factor $Y_O$ out from $G_{E(\mathcal D)}$ we are left
with an Euler product $A_{E(\mathcal D)}(\alpha,\gamma)$
which is absolutely convergent for
all of the variables in small disks around 0. Specifically,
\begin{align}
A_{E(\mathcal D)}(\alpha,\gamma)
 =\prod_p &\frac{\prod_{1\le j<k\le K}  (1-1/p^{1+\alpha_j+\alpha_k})
\prod_{1\le q< r\le Q}(1-1/p^
{1+\gamma_q+\gamma_r})
\prod_{q=1}^Q(1-1/p^{1+2\gamma_q})
}{\prod_{k=1}^K\prod_{q=1}^Q(1-1/p^{1+\alpha_k+\gamma_q})
  }\cr
& \qquad \times
\bigg(1+(1+\tfrac 1p)^{-1}
\sum_{0< \sum_k a_k +\sum_q c_q  \mbox{ is even}}
\frac{\prod_k \lambda(p^{a_k})\prod_q \mu_E(p^{c_q})}
{  p ^{\sum_k a_k(1/2+\alpha_k) +\sum_q c_q(1/2+\gamma_q)}}\bigg).
\end{align}
Note: when the above expression is evaluated, primes dividing $M$ contribute differently
 than primes not dividing $M$. These cases are analyzed when we refine our expression for $A$
in section~\ref{sec:AED}.

The functional equation may be written as
\begin{equation}
L_{E}(s,\chi_d)=w(E_d)
\bigg(\frac{\sqrt{M}|d|}{2\pi}\bigg)^{1-2s}g (s) L_{E}(1-s,\chi_d)
\end{equation}
where
\begin{equation}
g (s)=\frac{\Gamma(1-s)}{\Gamma(s)}.
\end{equation}
Note that $g (1/2)=1$.
The analytic conductor of $L_E(s,\chi_d)$ is $\frac{M|d|^2}{4 \pi^2}$ so that
the role of $2N$ is Theorem~\ref{thm:Sp} is played by
$2N=\log \frac{M|d|^2}{4 \pi^2}$.
(Again there are some subtleties concerning the ``conductor''.
See the discussion of Conjecture~1.5.5 in~\cite{CFKRS2}).
We are led to
\begin{conjecture}\label{conj:newp}
Suppose that the real parts of $\alpha_k$ and
$\gamma_q$ are positive. Then
\begin{align}
\sum_{|d|\le X\atop w(E_d)=1} &\frac{\prod_{k=1}^K L_E(1/2+\alpha_k,\chi_d)}
{\prod_{q=1}^Q L_E(1/2+\gamma_q,\chi_d)}\cr
&\qquad  =
\sum_{\epsilon \in \{-1,1\}^K}
Y_O A_{E(\mathcal D)}(\epsilon_1 \alpha_{1},\dots, \epsilon_K  \alpha_{K};
 \gamma )\prod_{k=1}^K g \big(\frac 12+\frac{\alpha_k-\epsilon_k \alpha_k }{2}\big)
\cr
 &\qquad \qquad \times
\sum_{|d|\le X \atop w(E_d)=1}
\bigg(\frac{M|d|^2}{4 \pi^2}\bigg)^{\frac 1 2\sum_{k=1}^K(\epsilon_k\alpha_k-\alpha_k)}
 +O(X^{1/2+\epsilon}).
\end{align}
\end{conjecture}
A conjecture for the moments of ratios of $L_E(s,\chi_d)$ over those $d$ for which
$w(E_d)=-1$  can be
analogously formulated:
\begin{conjecture}\label{conj:newm} Suppose that the real parts of $\alpha_k$ and
$\gamma_q$ are positive. Then
\begin{align}
\sum_{|d|\le X\atop w(E_d)=-1} &\frac{\prod_{k=1}^K L_E(1/2+\alpha_k,\chi_d)}
{\prod_{q=1}^Q L_E(1/2+\gamma_q,\chi_d)}\cr
&\qquad  =
\sum_{\epsilon \in \{-1,1\}^K} \mbox{sgn}(\epsilon)
Y_O A_{E(\mathcal D)}(\epsilon_1 \alpha_{1},\dots, \epsilon_K  \alpha_{K};
 \gamma )\prod_{k=1}^K g \big(\frac 12+\frac{\epsilon_k \alpha_k-\alpha_k}{2}\big)
\cr
 &\qquad \qquad \times
\sum_{|d|\le X \atop w(E_d)=-1}
\bigg(\frac{M|d|^2}{4 \pi^2}\bigg)^{\frac 12\sum_{k=1}^K(\epsilon_k\alpha_k-\alpha_k)}
 +O(X^{1/2+\epsilon}).
\end{align}
\end{conjecture}

\section{Refinements of the conjectures}\label{sec:refinements}
In this section we refine our conjectures in two ways.
We find closed form expressions
for the Euler products $A$ and we
express the combinatorial sums in our conjectures as
residues of multiple integrals.
This is similar to the treatment in~\cite{CFKRS2}.

\subsection{Closed form expressions for $A_\zeta$}\label{sec:closedformzeta}
Let $e(\theta)=e^{2 \pi i \theta}$ and let
$\delta_0(n)$ be the function which is 1 when $n=0$ and is 0 otherwise.
Then  $\delta_0(n)=\int_0^1e(n\theta)~d\theta .$
In the formula for $A_\zeta$ we then replace the summation condition
 $\sum a_k+\sum c_q=\sum b_\ell+\sum d_r$
by
\begin{align}
&\delta_0\big(\sum a_k+\sum c_q-\sum b_\ell-\sum d_r\big) \cr
&\phantom{XXXXXXXXXXX}=\int_0^1 e\bigg(\big(\sum a_k+\sum c_q-\sum b_\ell-\sum d_r\big)\theta\bigg)~d\theta.
\end{align}
After summing the geometric series that arise we deduce
\begin{lemma} Let $e(\theta)=e^{2 \pi i \theta}$. Then
\begin{align}
\sum_{\sum  a_k+\sum c_q=\sum b_\ell+
\sum d_r}&
\frac{\prod  \mu(p^c_q)\prod \mu(p^{d_r})}{  p^{\sum(1/2+\alpha_k)
a_k+\sum (1/2+\beta_\ell)b_\ell+\sum(1/2+\gamma_q)c_q+\sum(1/2+\delta_r)d_r}}\cr
&\phantom{XXXXXXXX} =
\int_0^1\frac{ \prod_{q=1}^Q\big(1-\frac{e(\theta)}{p^{1/2+\gamma_q}}\big)
\prod_{r=1}^R\big(1-\frac{e(-\theta)}{p^{1/2+\delta_r}}\big)}
{\prod_{k=1}^K\big(1-\frac{e(\theta)}{p^{1/2+\alpha_k}}\big)
\prod_{\ell=1}^L\big(1-\frac{e(-\theta)}{p^{1/2+\beta_\ell}}\big)}~d\theta .
\end{align}
\end{lemma}

\begin{corollary}
\begin{align}
A_\zeta(\alpha;\beta;\gamma;\delta)=&\prod_p \frac{\prod_{K=1 }^K\prod_{\ell=1}^L (1-1/p^{1+\alpha_k+\beta_\ell})
\prod_{q=1}^Q\prod_{r=1}^R(1-1/p^
{1+\gamma_q+\delta_r})}{\prod_{k=1}^K\prod_{r=1}^R(1-1/p^{1+\alpha_k+\delta_r})
\prod_{\ell=1}^L\prod_{q=1}^Q(1-1/p^{1+\beta_\ell+\gamma_q}) }\cr
& \qquad \times
\int_0^1\frac{ \prod_{q=1}^Q\big(1-\frac{e(\theta)}{p^{1/2+\gamma_q}}\big)
\prod_{r=1}^R\big(1-\frac{e(-\theta)}{p^{1/2+\delta_r}}\big)}
{\prod_{k=1}^K\big(1-\frac{e(\theta)}{p^{1/2+\alpha_k}}\big)
\prod_{\ell=1}^L\big(1-\frac{e(-\theta)}{p^{1/2+\beta_\ell}}\big)}~d\theta .
\end{align}
\end{corollary}

\subsection{Closed form expressions for $A_\mathcal D$}
Suppose that $f(x)=1+\sum_{n=1}^\infty u_n x^n$. Then
 \begin{equation}
 \sum_{0< n \mbox{ is even}} u_nx^n =\frac 12 \big(f(x)+f(-x)-2  \big)
\end{equation}
and
\begin{equation}
 1+(1+\tfrac 1p)^{-1}\sum_{0< n \mbox{ is even}} u_nx^n
=
\frac{1}{1+\tfrac 1 p}\bigg(\frac{f(x)+f(-x)}{2}+\frac 1 p \bigg).
\end{equation}
We apply this with
\begin{equation}
f(1/p)=\sum_{  a_k  , c_q  }
\frac{\prod_q \mu(p^{c_q})}
{  p ^{\sum_k a_k(1/2+\alpha_k) +\sum_q c_q(1/2+\gamma_q)}}
=\frac{\prod_{q=1}^Q\big(1-\frac{1}{p^{1/2+\gamma_q}}\big)}
{\prod_{k=1}^K\big(1-\frac{1}{p^{1/2+\alpha_k}}\big)}
\end{equation}
to deduce
\begin{lemma}
\begin{align}
 1+(1+\tfrac 1p)^{-1}&
\sum_{0< \sum_k a_k +\sum_q c_q  \mbox{ is even}}
\frac{\prod_q \mu(p^{c_q})}
{  p ^{\sum_k a_k(1/2+\alpha_k) +\sum_q c_q(1/2+\gamma_q)}} \cr
&\qquad =
\frac{1}{1+\tfrac 1 p}\bigg(\frac 1 2
\frac{\prod_{q=1}^Q\big(1-\frac{1}{p^{1/2+\gamma_q}}\big)}
{\prod_{k=1}^K\big(1-\frac{1}{p^{1/2+\alpha_k}}\big)}
+\frac 1 2
\frac{\prod_{q=1}^Q\big(1+\frac{1}{p^{1/2+\gamma_q}}\big)}
{\prod_{k=1}^K\big(1+\frac{1}{p^{1/2+\alpha_k}}\big)}
+\frac 1 p \bigg).
\end{align}
\end{lemma}

\begin{corollary}
\begin{align}
A_\mathcal D(\alpha,\gamma)
 =&\prod_p \frac{\prod_{  j\le k\le K}  (1-1/p^{1+\alpha_j+\alpha_k})
\prod_{ q< r\le Q}(1-1/p^
{1+\gamma_q+\gamma_r})}{\prod_{k=1}^K\prod_{q=1}^Q(1-1/p^{1+\alpha_k+\gamma_q})
  }\cr
& \qquad \times
\frac{1}{1+\tfrac 1 p}\bigg(\frac 1 2
\frac{\prod_{q=1}^Q\big(1-\frac{1}{p^{1/2+\gamma_q}}\big)}
{\prod_{k=1}^K\big(1-\frac{1}{p^{1/2+\alpha_k}}\big)}
+\frac 1 2
\frac{\prod_{q=1}^Q\big(1+\frac{1}{p^{1/2+\gamma_q}}\big)}
{\prod_{k=1}^K\big(1+\frac{1}{p^{1/2+\alpha_k}}\big)}
+\frac 1 p \bigg).
\end{align}
\end{corollary}

\subsection{Closed form expressions for $A_{E(\mathcal D)}$}\label{sec:AED}
For simplicity, let $E=E_{11}$.
We apply the method of the last section, this time with
\begin{equation}
f(1/p)=\sum_{  a_k  , c_q  }
\frac{\prod_k \lambda(p^{a_k})\prod_q \mu_E(p^{c_q})}
{  p ^{\sum_k a_k(1/2+\alpha_k) +\sum_q c_q(1/2+\gamma_q)}} .
\end{equation}
If $p$ is not 11, then
\begin{equation}
f(1/p)=\frac{\prod_{q=1}^Q\big(1-\frac{\lambda(p)}{p^{1/2+\gamma_q}}
+\frac{1}{p^{1+2\gamma_q}}\big)}
{\prod_{k=1}^K\big(1-\frac{\lambda(p)}{p^{1/2+\alpha_k}}+
\frac{\lambda(p)}{p^{1+2\alpha_k}}\big)} ,
\end{equation}
whereas if $p=11$, then
\begin{equation}
f(1/p)=\frac{\prod_{q=1}^Q\big(1-\frac{1/\sqrt{11}}{11^{1/2+\gamma_q}}
 \big)}
{\prod_{k=1}^K\big(1-\frac{1/\sqrt{11}}{11^{1/2+\alpha_k}}
 \big)}
\end{equation}

\begin{lemma}If $p\ne 11$, then
\begin{align}
 1+(1+\tfrac 1p)^{-1}&
\sum_{0< \sum_k a_k +\sum_q c_q  \mbox{ is even}}
\frac{\prod_k \lambda(p^{a_k})\prod_q \mu_E(p^{c_q})}
{  p ^{\sum_k a_k(1/2+\alpha_k) +\sum_q c_q(1/2+\gamma_q)}}
 \cr
&\qquad =
\frac{1}{1+\tfrac 1 p}\bigg(\frac 1 2
\frac{\prod_{q=1}^Q\big(1-\frac{\lambda(p)}{p^{1/2+\gamma_q}}
+\frac{1}{p^{1+2\gamma_q}}\big)}
{\prod_{k=1}^K\big(1-\frac{1}{p^{1/2+\alpha_k}}+
\frac{\lambda(p)}{p^{1+2\alpha_k}}\big)}+
\frac 1 2
\frac{\prod_{q=1}^Q\big(1+\frac{\lambda(p)}{p^{1/2+\gamma_q}}
+\frac{1}{p^{1+2\gamma_q}}\big)}
{\prod_{k=1}^K\big(1+\frac{\lambda(p)}{p^{1/2+\alpha_k}}+
\frac{\lambda(p)}{p^{1+2\alpha_k}}\big)}
+\frac 1 p \bigg).
\end{align}
while if $p=11$, then
\begin{align}
 1+(1+\tfrac 1p)^{-1}&
\sum_{0< \sum_k a_k +\sum_q c_q  \mbox{ is even}}
\frac{\prod_k \lambda(p^{a_k})\prod_q \mu_E(p^{c_q})}
{  p ^{\sum_k a_k(1/2+\alpha_k) +\sum_q c_q(1/2+\gamma_q)}}
 \cr
&\qquad =
\frac{1}{1+\tfrac 1 {11}}\bigg(\frac 1 2
\frac{\prod_{q=1}^Q\big(1-\frac{1/\sqrt{11}}{11^{1/2+\gamma_q}}
 \big)}
{\prod_{k=1}^K\big(1-\frac{1/\sqrt{11}}{11^{1/2+\alpha_k}}
 \big)}+
\frac 1 2
\frac{\prod_{q=1}^Q\big(1+\frac{1/\sqrt{11}}{11^{1/2+\gamma_q}}
 \big)}
{\prod_{k=1}^K\big(1+\frac{1/\sqrt{11}}{11^{1/2+\alpha_k}}
 \big)}
+\frac 1 {11} \bigg).
\end{align}

\end{lemma}

\begin{corollary}
\begin{align}
A_{E_{11}(\mathcal D)}&(\alpha,\gamma)\cr
& =
 \frac{\prod_{1\le j<k\le K}  (1-1/11^{1+\alpha_j+\alpha_k})
\prod_{1\le q< r\le Q}(1-1/11^
{1+\gamma_q+\gamma_r})}{\prod_{k=1}^K\prod_{q=1}^Q(1-1/11^{1+\alpha_k+\gamma_q})
\prod_{q=1}^Q(1-1/11^{1+2\gamma_q})
  }\cr
&   \quad\times
\frac{1}{1+\tfrac 1 {11}}\bigg(\frac 1 2
\frac{\prod_{q=1}^Q\big(1-\frac{1/\sqrt{11}}{11^{1/2+\gamma_q}}
 \big)}
{\prod_{k=1}^K\big(1-\frac{1/\sqrt{11}}{11^{1/2+\alpha_k}}
 \big)}+
\frac 1 2
\frac{\prod_{q=1}^Q\big(1+\frac{1/\sqrt{11}}{11^{1/2+\gamma_q}}
 \big)}
{\prod_{k=1}^K\big(1+\frac{1/\sqrt{11}}{11^{1/2+\alpha_k}}
 \big)}
+\frac 1 {11} \bigg)\cr
&  \quad\times
 \prod_{p\ne 11} \frac{\prod_{1\le j<k\le K}  (1-1/p^{1+\alpha_j+\alpha_k})
\prod_{1\le q< r\le Q}(1-1/p^
{1+\gamma_q+\gamma_r})}{\prod_{k=1}^K\prod_{q=1}^Q(1-1/p^{1+\alpha_k+\gamma_q})
\prod_{q=1}^Q(1-1/p^{1+2\gamma_q})
  }\cr
&   \qquad\times
\frac{1}{1+\tfrac 1 p}\bigg(\frac 1 2
\frac{\prod_{q=1}^Q\big(1-\frac{\lambda(p)}{p^{1/2+\gamma_q}}
+\frac{1}{p^{1+2\gamma_q}}\big)}
{\prod_{k=1}^K\big(1-\frac{1}{p^{1/2+\alpha_k}}+
\frac{\lambda(p)}{p^{1+2\alpha_k}}\big)}+
\frac 1 2
\frac{\prod_{q=1}^Q\big(1+\frac{\lambda(p)}{p^{1/2+\gamma_q}}
+\frac{1}{p^{1+2\gamma_q}}\big)}
{\prod_{k=1}^K\big(1+\frac{\lambda(p)}{p^{1/2+\alpha_k}}+
\frac{\lambda(p)}{p^{1+2\alpha_k}}\big)}
+\frac 1 p \bigg)
 .
\end{align}
\end{corollary}

\subsection{Combinatorial sums as integrals}\label{sec:combsum}

We express the sums appearing in our conjectures in terms of
multiple integrals.  The expressions will involve the
Vandermonde determinant,
given by
 \begin{equation}
\Delta(w_1,\dots,w_R)=\det_{R\times R}\big( w_i^{j-1}\big).
\end{equation}
We often omit the subscripts and write $\Delta(w)$ in place of
$\Delta(w_1,\dots,w_R)$.
The key fact about the Vandermonde is that
\begin{equation} \Delta(w_1,\dots,w_R)=\prod_{1\le j < k\le R}(w_j-w_i).
\end{equation}

\begin{lemma}
Suppose that
$F( z;w)=F( z_1,\dots, z_K;w_1,\dots,w_L)$
is a function of $K+L$ variables,
which is symmetric with respect to the first $K$ variables and
symmetric with respect to the second set of $L$ variables.
Suppose also that $F$ is regular near $(0,\dots,0)$.
Suppose further that $f(s)$ has a simple pole
of residue~$1$ at $s=0$
but is otherwise analytic in $|s|\le 1$.
Let
\begin{equation}
H( z_1,\dots, z_K;w_1,\dots w_L)=
F( z_1,\dots;\dots,w_L)
\prod_{k=1}^K\prod_{\ell=1}^L
f( z_k-w_\ell).\end{equation}
If $| \alpha_k|<1$ then
\begin{align}
\sum_{\sigma \in \Xi_{K,L}} &
H( \alpha_{\sigma(1)},\dots,  \alpha_{\sigma(K)};
\alpha_{\sigma(K+1)}\dots  \alpha_{\sigma(K+L)})=\cr
&\qquad
\frac{(-1)^{(K+L)(K+L-1)/2}}{K!L! (2\pi i)^{K+L}}
\int\limits_{|z_i|=1}\frac{H(z_1, \dots,z_K;z_{K+1},\dots,z_{K+L})
\Delta(z_1,\dots,z_{K+L})^2}{\prod_{j=1}^{K+L}\prod_{k=1}^{K+L}
(z_k-\alpha_j)}
\,dz_1\dots dz_{K+L}.
\end{align}
\end{lemma}

In view of the last formula of Section~\ref{sec:closedformzeta},
and using
\begin{equation}
H_{\zeta,t}(w_1, \dots,w_{K+L};\gamma;\delta)=
\frac{\prod_{\ell=1}^L\chi(s-w_{K+\ell})^{1/2}}{\prod_{k=1}^K \chi(s+w_k)^{1/2}}
Y_UA_\zeta(w_1,\dots,w_K; -w_{K+1},\dots ,-w_{K+L};\gamma;\delta)
\end{equation}
with
\begin{align}
A_\zeta(\alpha;\beta;\gamma;\delta)=&\prod_p \frac{\prod_{K=1 }^K\prod_{\ell=1}^L (1-1/p^{1+\alpha_k+\beta_\ell})
\prod_{q=1}^Q\prod_{r=1}^R(1-1/p^
{1+\gamma_q+\delta_r})}{\prod_{k=1}^K\prod_{r=1}^R(1-1/p^{1+\alpha_k+\delta_r})
\prod_{\ell=1}^L\prod_{q=1}^Q(1-1/p^{1+\beta_\ell+\gamma_q}) }\cr
& \qquad \times
\int_0^1\frac{ \prod_{q=1}^Q\big(1-\frac{e(\theta)}{p^{1/2+\gamma_q}}\big)
\prod_{r=1}^R\big(1-\frac{e(-\theta)}{p^{1/2+\delta_r}}\big)}
{\prod_{k=1}^K\big(1-\frac{e(\theta)}{p^{1/2+\alpha_k}}\big)
\prod_{\ell=1}^L\big(1-\frac{e(-\theta)}{p^{1/2+\beta_\ell}}\big)}~d\theta .
\end{align}
and
\begin{equation}Y_U(\alpha;\beta;\gamma;\delta):=
\frac{\prod_{k=1}^K\prod_{\ell=1}^L \zeta(1+\alpha_k+\beta_\ell)\prod_{q=1}^Q
\prod_{r=1}^R\zeta(1+\gamma_q+\delta_r)}
{\prod_{k=1}^K\prod_{r=1}^R \zeta(1+\alpha_k+\delta_r)
\prod_{\ell=1}^L\prod_{q=1}^Q \zeta(1+\beta_\ell+\gamma_q)},\end{equation}
we can reformulate
Conjecture~\ref{conj:zeta} as
 \begin{align}
\int_0^T &
\frac{\prod_{k=1}^{K}\zeta(s+\alpha_k)\prod_{\ell=K+1}^{K+L}
\zeta(1-s-\alpha_\ell)}{\prod_{q=1}^Q\zeta(s+\gamma_q)
\prod_{r=1}^R\zeta(1-s+\delta_r)} ~dt = \int_0^T
\frac{\prod_{k=1}^K\chi(s+\alpha_k)^{1/2}}{\prod_{\ell=1}^L \chi(s-\alpha_{K+\ell})^{1/2}}
\cr
&\qquad
\times
\frac{(-1)^{(K+L)(K+L-1)/2}}{K!L! (2\pi i)^{K+L}}
\int\limits_{|z_i|=1}\frac{H_{\zeta,t}(z_1, \dots,z_K;z_{K+1},\dots,z_{K+L})
\Delta(z_1,\dots,z_{K+L})^2}{\prod_{j=1}^{K+L}\prod_{k=1}^{K+L}
(z_k-\alpha_j)}
\,dz_1\dots dz_{K+L}
 ~dt\cr
& \qquad \qquad
+O(T^{1/2+\epsilon}).
\end{align}
This should be compared with the reformulation of Theorem~\ref{thm:U}:
\begin{align}
\int_{U(N)}&\frac{\prod_{k=1}^K\Lambda_A(e^{-\alpha_j})\prod_{\ell=K+1}^{K+L}
\Lambda_{A^*}(e^{\alpha_\ell})} {\prod_{q=1}^Q\Lambda_A(e^{-\gamma_q})
\prod_{r= 1}^{R}
\Lambda_{A^*}(e^{-\delta_r})}dA_N
= e^{\frac N2(- \sum_{k=1}^K   \alpha_k+\sum_{\ell=1}^L
\alpha_{K+\ell})}
\cr
& \times\frac{(-1)^{(K+L)(K+L-1)/2}}{K!L! (2\pi i)^{K+L}}
\int\limits_{|z_i|=1}
 \frac{H_U(z_1, \dots,z_K;z_{K+1},\dots,z_{K+L};\gamma;\delta)
\Delta(z_1,\dots,z_{K+L})^2}{\prod_{j=1}^{K+L}\prod_{k=1}^{K+L}
(z_k-\alpha_j)}
\prod_k dz_k
\end{align}
where
\begin{equation}
H_U(w_1, \dots,w_{K+L};\gamma;\delta)=
e^{\frac N2 \sum_{k=1}^K   w_k -\frac N 2
\sum_{\ell=1}^L w_{K+\ell} }
y_U(w_1,\dots,w_K; -w_{K+1},\dots ,-w_{K+L};\gamma;\delta).
\end{equation}

Note that for a small shift $\alpha$,
\begin{eqnarray}
\chi(s+\alpha)=\bigg(\frac{|t|}{2\pi}\bigg)^{1/2-s-\alpha}\big(1+O(1/(1+|t|)\big)
\end{eqnarray}
so that, for example,
\begin{eqnarray*}
\prod_{k=1}^K\frac{\chi(s+\alpha_k)^{1/2}}{\chi(s)^{1/2}}\prod_{\ell=1}^L \frac{\chi(s-\alpha_{K+\ell})^{1/2}}{\chi(s)^{1/2}}=e^{\frac \ell 2(- \sum_{k=1}^K   \alpha_k+\sum_{\ell=1}^L
\alpha_{K+\ell})}\big(1+O(1/(1+|t|)\big)
\end{eqnarray*}
where $\ell=\log \frac {t}{2\pi}$, which compares with the random matrix formula with $N$ replaced by $\ell$.
For large shifts, this approximation deteriorates and necessitates that we retain the accurate expression
involving the product of $\chi$.

Note also that this integral formula gives an analytic continuation in the
variables $\alpha$ and $\beta$ (so that we no longer need to restrict them to have
positive real parts), whereas the variables $\gamma$ and $\delta$ are still
required to have positive real parts.
\begin{lemma}
Suppose that
$F(z )=F(z_1,\dots,z_K )$
is a function of $K $ variables,
which is symmetric and regular near $(0,\dots,0)$.
Suppose further that $f(s)$ has a simple pole
of residue~$1$ at $s=0$
but is otherwise analytic in $|s|\le 1$.
Let
either
\begin{equation}
H(z_1,\dots,z_K )=
F(z_1,\dots,z_K)
\prod_{1\le j\le k\le K}
f(z_j+z_k)\end{equation}
or
\begin{equation}
H(z_1,\dots,z_K )=
F(z_1,\dots,z_K)
\prod_{1\le j< k\le K}
f(z_j+z_k).\end{equation}
If $|\alpha_k|<1$ then
\begin{align}
\sum_{\epsilon \in \{-1,+1\}^K}&
   H(\epsilon_1 \alpha_1,\dots,\epsilon_K \alpha_K)\cr
&\qquad =
\frac{(-1)^{K(K-1)/2}2^K}{K!  (2\pi i)^{K }}
\int\limits_{|z_i|=1}\frac{H(z_1, \dots,z_K )
\Delta(z_1^2,\dots,z_{K }^2)^2\prod_{k=1}^K z_k}
{\prod_{j=1}^{K }\prod_{k=1}^{K }
(z_k-\alpha_j)(z_k+\alpha_j)}
\,dz_1\dots dz_{K }
\end{align}
and
\begin{align}
\sum_{\epsilon \in \{-1,+1\}^K}&\mbox{sgn}(\epsilon)
H(\epsilon_1 \alpha_1,\dots,\epsilon_K \alpha_K)\cr
&\qquad =
\frac{(-1)^{K(K-1)/2}2^K}{K!  (2\pi i)^{K }}
\int\limits_{|z_i|=1}\frac{H(z_1, \dots,z_K )
\Delta(z_1^2,\dots,z_{K }^2)^2\prod_{k=1}^K \alpha_k}
{\prod_{j=1}^{K }\prod_{k=1}^{K }
(z_k-\alpha_j)(z_k+\alpha_j)}
\,dz_1\dots dz_{K }.
\end{align}
\end{lemma}

Using this Lemma, we can reformulate Theorem~\ref{thm:Sp} and
Conjecture~\ref{conj:realchi} as
\begin{align}
\int_{USp(2N)}&\frac{\prod_{k=1}^K\Lambda_A(e^{-\alpha_k})}
{\prod_{q=1}^Q \Lambda_A(e^{-\gamma_q})}~dA
=e^{-\frac N 2 \sum_{k=1}^K \alpha_k}
\frac{(-1)^{K(K-1)/2}2^K}{K!  (2\pi i)^{K }} \cr
& \times \int\limits_{|z_i|=1}\frac{h_S(z_1, \dots,z_K;\gamma)
\Delta(z_1^2,\dots,z_{K }^2)^2\prod_{k=1}^K z_k}
{\prod_{j=1}^{K }\prod_{k=1}^{K }
(z_k-\alpha_j)(z_k+\alpha_j)}
\,dz_1\dots dz_{K },
\end{align}
and
\begin{align}
\sum_{0<d\le X} &\frac{\prod_{k=1}^K L(1/2+\alpha_k,\chi_d)}
{\prod_{q=1}^Q L(1/2+\gamma_q,\chi_d)}=\sum_{0<d\le X}
\bigg(\frac{|d|}{\pi}\bigg)^{-\frac 12\sum_{k=1}^K \alpha_k}\frac{(-1)^{K(K-1)/2}2^K}{K!
(2\pi i)^{K }}
\cr
&\qquad
\times \int\limits_{|z_i|=1}\frac{H_{\mathcal D+,d,\alpha,\gamma}(z_1, \dots,z_K;\gamma )
\Delta(z_1^2,\dots,z_{K }^2)^2\prod_{k=1}^K z_k}
{\prod_{j=1}^{K }\prod_{k=1}^{K }
(z_k-\alpha_j)(z_k+\alpha_j)}
\,dz_1\dots dz_{K }
 +O(X^{1/2+\epsilon}).
\end{align}

There is a similar reformulation of Theorem~\ref{thm:Oe} and~\ref{thm:Oo} and
Conjectures~\ref{conj:newp} and~\ref{conj:newm}.

\section{Examples and Applications}\label{sec:applications}

\subsection{Farmer's conjecture revisited}

We first give a more precise version of Farmer's original conjecture.
We use Conjecture~\ref{conj:zeta} with $K=L=1$. In this case,
$\Xi_{1,1}=\{(1),(12)\}$ consists of the identity permutation
and the transposition $(12)$. We identify $\alpha_1=\alpha$, $\alpha_2=-\beta$,
$\gamma_1=\gamma$, and $\delta_1=\delta$. Then Conjecture~\ref{conj:zeta} tells us that
\begin{align}\label{eqn:farmer}
 \int_0^T &\frac{\zeta(s+\alpha)\zeta(1-s+\beta)}
{\zeta(s+\gamma)\zeta(1-s+\delta)}~dt \\
&\qquad =\int_0^T \bigg(Y_UA_\zeta(\alpha,\beta;\gamma;\delta)
+\bigg(\frac{t}{2 \pi}\bigg)^{-\alpha -\beta} \nonumber
Y_UA_\zeta(-\beta, -\alpha;\gamma;\delta)\bigg)~dt+O(T^{1/2+\epsilon}) .
\end{align}
We see that
\begin{equation}
Y_U(\alpha;\beta;\gamma;\delta)=\frac{\zeta(1+\alpha+\beta)\zeta(1+\gamma+\delta)}
{\zeta(1+\alpha+\delta) \zeta(1+\beta+\gamma)}
\end{equation}
and
\begin{equation}
A_\zeta(\alpha;\beta;\gamma;\delta)
=\prod_p \frac{  (1-1/p^{1+\alpha +\beta })
(1-1/p^{1+\gamma +\delta })}
{ (1-1/p^{1+\alpha +\delta }) (1-1/p^{1+\beta +\gamma ) }}
\int_0^1\frac{  \big(1-\frac{e(\theta)}{p^{1/2+\gamma }}\big)
\big(1-\frac{e(-\theta)}{p^{1/2+\delta }}\big)}
{ \big(1-\frac{e(\theta)}{p^{1/2+\alpha }}\big)
\big(1-\frac{e(-\theta)}{p^{1/2+\beta }}\big)}~d\theta .
\end{equation}

For values of $\alpha, \beta, \gamma, \delta \to 0$ we
have asymptotically
\begin{equation}
Y_U(\alpha;\beta;\gamma;\delta)\sim \frac
{ ( \alpha+\delta)  ( \beta+\gamma)}
{ ( \alpha+\beta) ( \gamma+\delta)},
\end{equation}
$A_\zeta\sim 1$, and $\frac{t}{2\pi}$ can be replaced by $T$;
to a first order approximation, we then have
\begin{align}
\frac{1}{T}\int_0^T \frac{\zeta(s+\alpha)\zeta(1-s+\beta)}
{\zeta(s+\gamma)\zeta(1-s+\delta)}~dt
\sim& \frac{(\alpha+\delta)(\beta+\gamma)}{(\alpha+\beta)(\gamma+\delta)}
+T^{-\alpha-\beta}\frac{(-\beta+\delta)(-\alpha+\gamma)}{(-\beta-\alpha)
(\gamma+\delta)}\cr
=& 1+(1-T^{-\alpha-\beta})\frac{(\alpha-\gamma)(\beta-\delta)}{(\alpha+\beta)
(\gamma+\delta)},
\end{align}
which recovers Farmer's original conjecture.

\subsection{Logarithmic derivatives of $\zeta(s)$}
 Goldston,
Gonek, and Montgomery [GGM] proved, assuming the Riemann Hypothesis, that
\begin{equation}
\frac{1}{T}
\int_0^T \bigg| \frac{\zeta'(1/2+r+it)}{\zeta(1/2+r+it)}\bigg|^2 ~dt
\sim \sum_p \frac{\log ^2p}{p^{1+2r}-1}-\frac{T^{-2r}}{4r^2}
+ \log^2T \int_1^\infty (F(\alpha,T)-1)T^{-2 r\alpha} ~d\alpha
\end{equation}
uniformly for
$T^{-1/2}\log T \ll r \ll 1 $, where
$F(\alpha,T)$ is Montgomery's pair correlation function.
Montgomery's function is expected to
satisfy $F(\alpha,T)=1+o(1)$ uniformly for
bounded $\alpha$ so that the term involving $F$ is expected
to be small.  Also, the sum over primes is $ ~\sim \frac {1}{4r^2}$
as $r\to 0$.

Here we obtain a conjecture for this quantity
which is more precise than the Goldston-Gonek-Montgomery
formula, in that it contains some
lower order terms and we expect it to
be accurate with a square-root error term.
We deduce our conjecture by differentiating
the formula of the last section with respect to $\alpha$ and $\beta$ and setting
$\alpha=\beta=\gamma=\delta=r$.  To help compute this, the following formulas,
about a function $f$ which is analytic in a neighborhood
of the origin,
are helpful:
\begin{equation}
\frac{d}{da}  \frac{f(a+b)f(c+d)}{f(a+d)f(b+c)}\bigg|_{a=b=c=d=r}
=\frac{d}{db}  \frac{f(a+b)f(c+d)}{f(a+d)f(b+c)}\bigg|_{a=b=c=d=r}
=0,
\end{equation}
and
\begin{align}
&\frac{d}{da} \frac{d}{db} \frac{f(a+b)f(c+d)}{f(a+d)f(b+c)}\bigg|_{a=b=c=d=r}\cr
&\qquad
=\frac{f''(2r)}{f(2r)}-\bigg(\frac{f'(2r)}{f(2r)}\bigg)^2
=\frac d{dx} \frac{f'(x)}{f(x)}\bigg|_{x=2r}=\frac{d^2}{dx^2}\log(f(x))\bigg|_{x=2r},
\end{align}
Also,
 \begin{equation}
\frac {d}{da} \frac d {db} \frac{f(a)f(b)}{f(c)f(d)}\bigg|_{a=b=c=d=r}
=\bigg(\frac{f'}{f}(2r)\bigg)^2.
\end{equation}
Thus, we can now calculate
\begin{align}
&\frac{d}{d\alpha}\frac {d}{d\beta}Y_UA_\zeta(\alpha;\beta;\gamma;\delta)
\bigg|_{\alpha=\beta=\gamma=\delta=r}\cr
&\qquad
=\bigg(\frac{\zeta'}{\zeta}\bigg)'(1+2r) +\sum_p
\bigg( \frac{-p^{1+2r}\log^2 p }{(p^{1+2r}-1)^2}+\int_0^1\frac{\log ^2p}
{(e(\theta)p^{\frac 12 +r} -1)^2}~d\theta\bigg)\ ,
\end{align}
and
\begin{align}
&\frac{d}{d\alpha}\frac {d}{d\beta}\bigg(\frac{t}{2\pi}\bigg)^{-\alpha-\beta}
Y_UA_\zeta(-\beta;-\alpha;\gamma;\delta)
\bigg|_{\alpha=\beta=\gamma=\delta=r}\cr
&\qquad
= \bigg(\frac{t}{2\pi}\bigg)^{-2r}
A_\zeta(-r,-r,r,r)\zeta(1-2r)\zeta(1+2r).
\end{align}

Thus, we have
\begin{conjecture}
\begin{align}
\frac{1}{T}
\int_0^T &\bigg| \frac{\zeta'(1/2+r+it)}{\zeta(1/2+r+it)}\bigg|^2 ~dt
\cr
&\qquad = \bigg(\frac{\zeta'}{\zeta}\bigg)'(1+2r)+
 \bigg( \frac{T}{2\pi}\bigg)^{-2r}
A_\zeta(-r,-r,r,r)\frac{\zeta(1-2r)\zeta(1+2r)}{1-2r} \cr
&\qquad \qquad +c(r)+O(T^{-1/2+\epsilon}),
\end{align}
where $c(r)$ is  a function of $r$ which is uniformly bounded for $|r|<1/4-\epsilon$
and is given by
\begin{equation}
c(r)=
\sum_p
\bigg( \frac{-p^{1+2r}\log^2 p }{(p^{1+2r}-1)^2}+\int_0^1\frac{\log ^2p}
{(e(\theta)p^{\frac 12 +r} -1)^2}~d\theta\bigg).
\end{equation}
\end{conjecture}

\subsection{A conjecture of Keating and Snaith}

A conjecture of Keating and Snaith is
\begin{equation}\frac{1}{T}\int_0^T \frac{\zeta(\frac 12 +it)^K}{\zeta(\frac 12 -it)^K}
dt\sim G(1-K)G(1+K)b_K (\log T)^{-K^2}\end{equation}
where $G$ is Barnes double Gamma-function and
\begin{equation}b_K=\prod_p(1-1/p)^{-K^2}\sum_{j=0}^\infty \frac{\Gamma(1+K)
\Gamma(1-K)}{\Gamma(1+K-j)\Gamma(1-K-j)j!^2p^j}.\end{equation}

Note that if $K$ is a positive integer, then $b_K=0$. Here we consider
the case that $K$ is a positive integer, but integrate
a ratio of shifted zeta-functions.  For $\alpha_k$ and $\delta_k$
with positive real parts we have, as a consequence of Conjecture~\ref{conj:zeta},
\begin{equation}
\frac{1}{T}\int_0^T \prod_{k=1}^K\frac{\zeta(s+\alpha_k) }
{ \zeta(1-s+\delta_k)} ~dt =B(\alpha, \delta) \prod_{i,j=1}^K \zeta(1+\alpha_i+\delta_j)^{-1}
+O(T^{1/2+\epsilon})
\end{equation}
where
\begin{equation}
B(\alpha, \delta)=
\prod_p
\prod_{i,j=1}^K\left(1-\frac{1}{p^{1+\alpha_i+\delta_j}}\right)^{-1}
\sum_{\sum a_i=\sum d_i}\frac{\prod \mu(p^{d_j})}
{p^{\sum a_i(\frac 12+\alpha_i)+\sum d_i(\frac 12 +\delta_i)}}
.
\end{equation}
Note that the product over primes in $B$ is absolutely
convergent for sufficiently small values of
the shifts $\alpha_i,\delta_i$ and is equal to the $b_K$ in Keating
 and Snaith's formula when all of the shifts are 0.
Also, the size of this expression is about $(\log T)^{-K^2}$ when all
of the shifts have order of magnitude $1/\log T$.

\subsection{Discrete moments of $\zeta$}
Let $\rho=\beta+i\gamma$ stand for a typical complex zero
of the Riemann zeta-function. The Riemann Hypothesis, which we assume here,
asserts that all $\beta=1/2$. The number of $\gamma\le T$ was proven by Riemann
and von~Mangoldt
to equal $\frac{T}{2\pi} \log \frac{T}{2 \pi e}+O(\log T)$.
Chris Hughes has conjectured a formula for the leading
term of
\begin{equation}\sum_{\gamma \le T}|\zeta'(\rho)|^{s}\end{equation}
for complex $s$ with $\sigma>-3$. His conjecture is based on
an exact formula he proved for the analogous random matrix moment:
\begin{equation}
\int_{U(N)}\sum_{n=1}^N |\Lambda_A'(e^{i\theta_n})|^s dA_N
= \frac{G^2(\tfrac s 2 +2)G(N+s+2)G(N)}{G(s+3)G^2(N+\tfrac s 2 +1)}
\end{equation}
for $\Re s>-3$.  Here the $e^{i\theta_n}$ are the zeros
of the characteristic polynomial $\Lambda_A$ and
$G$ is the Barnes double gamma-function.  It should be possible to use our ratios conjecture to determine the lower order terms
of the discrete moments $\sum_{\gamma \le T}|\zeta'(\rho)|^{2k}$ for positive integer $k$.  We will return to this in a later paper.

Now, we compute a conjecture
for
\begin{equation}
D(a,c)=\sum_{\gamma\le T}
\frac {\zeta(\rho+a)}{\zeta(\rho+c)},
\end{equation}
where $\Re a, \Re c >0$, a sum which
was considered in Farmer's
paper [F1]. Farmer's conjecture for this sum is
\begin{equation} \label{eqn:farmerconj}
D(a,c)\sim \frac{T}{2\pi}\left(\log T +(1-T^{-a})\left(\frac 1 c -\frac 1 a\right)\right).
\end{equation}
We now give a more precise conjecture based on our ratios conjecture.

By Cauchy's formula,
\begin{equation}
D(a,c)=\frac{1}{2\pi i}\int_{\mathcal C} \frac{\zeta'}{\zeta}(s)\frac{\zeta(s+a)}{\zeta(s+c)}~ds
\end{equation}
where $\mathcal C$ is a tall, thin rectangular path (with vertices $1/2\pm \alpha$, $1/2\pm \alpha +iT$, which encloses the zeros $1/2+i\gamma$ for $0< \gamma < T$. Here $\alpha>0$ is smaller than the real parts of $a$ and $b$.
The parameter $T$ can be
slightly adjusted if necessary to conclude that the integrals on the horizontal portions of the path are $\ll T^\epsilon.$ Also, the integral on the right hand path $1/2+\alpha+it, 0\le t\le T$ is $\ll T^\epsilon$ as can be seen by moving the path of integration to the right of $\sigma=1$ and integrating term-by-term. Thus, we have
\begin{equation}
D(a,c)=\frac{-1}{2\pi }\int_0^T \frac{\zeta'}{\zeta}(1/2-\alpha+it)\frac{\zeta(1/2-\alpha+a+it)}{\zeta(1/2-\alpha+c+it)}~dt
+O(T^\epsilon).
\end{equation}
Now we use the functional equation
\begin{equation}
\frac{\zeta'}{\zeta}(s)=\frac{\chi'}{\chi}(s)-\frac{\zeta'}{\zeta}(1-s)
\end{equation}
and split the integral into two pieces. The part with the $\chi'/\chi$ can be treated much as the first integral (on the $1/2+\alpha$-line) by moving the path of integration to the right and into the region of absolute convergence of the Dirichlet series where we can integrate term-by-term.  Note that this Dirichlet series begins with a 1 whereas in the first integral the series had no constant term.  Thus, the contribution from this integral is
\begin{align}
=&\frac{-1}{2\pi} \int_0^T \frac{\chi'}{\chi}(2+it)~dt+O(T^\epsilon)\cr
=&\frac{1}{2\pi}\int_0^T \log\frac{t}{2\pi} ~dt +O(T^\epsilon).
\end{align}

Thus, we have
\begin{align}
D(a,b)=&\frac{1}{2\pi }\int_0^T \left(\log\frac{t}{2\pi}+ \frac{\zeta'}{\zeta}(1/2+\alpha-it)\frac{\zeta(1/2-\alpha+a+it)}{\zeta(1/2-\alpha+c+it)}\right)
~dt
+O(T^\epsilon)\cr
=&  \frac{1}{2\pi }\int_0^T  \left(\log\frac{t}{2\pi}+ \frac{d}{d \beta} \frac{\zeta(1-s+\alpha+\beta ) \zeta(s-\alpha+a )}{\zeta(1-s+\alpha)\zeta(s-\alpha+c )}\right)~dt
\bigg|_{\beta=0}
+O(T^\epsilon)\end{align}
where $s=1/2+it$.
By (\ref{eqn:farmer}), we have
\begin{align}
 \int_0^T   \frac{\zeta(1-s+\alpha+\beta ) \zeta(s-\alpha+a )}{\zeta(1-s+\alpha)\zeta(s-\alpha+c )} ~dt  =\mathstrut &
O(T^{1/2+\epsilon})+ \cr
 +\int_0^T \bigg(Y_UA_\zeta(a-\alpha,\alpha+\beta;c-\alpha,\alpha)
& +\left(\frac {t}{2\pi}\right)^{-a-\beta}Y_UA_\zeta(-\alpha-\beta,\alpha-a;c-\alpha,\alpha)\bigg)~dt  .
\end{align}
Now it can be easily calculated that
\begin{equation}
Y_U(x,y;z,w)=\frac{\zeta(1+x+y)\zeta(1+z+w)}{\zeta(1+x+w)\zeta(1+y+z)}
\end{equation}
and
\begin{equation}
A_\zeta(x,y;z,w)=
 \prod_p \frac{\left(1-\frac{1}{p^{1+z+w}}\right)
\left(1-\frac{1}{p^{1+y+z}}
-\frac{1}{p^{1+x+w}}+\frac{1}{p^{1+z+w}}\right)}
{\left(1-\frac{1}{p^{1+y+z}}\right)
\left(1-\frac{1}{p^{1+x+w}}\right)} .
\end{equation}
Thus,
\begin{equation} Y_UA_\zeta(a-\alpha,\alpha+\beta;c-\alpha,\alpha)
= \frac{\zeta(1+a+\beta)\zeta(1+c)}{\zeta(1+a)\zeta(1+c+\beta)}
\prod_p \frac{\left(1-\frac{1}{p^{1+c}}\right)
\left(1-\frac{1}{p^{1+\beta+c}}
-\frac{1}{p^{1+a}}+\frac{1}{p^{1+c}}\right)}
{\left(1-\frac{1}{p^{1+\beta+c}}\right)
\left(1-\frac{1}{p^{1+a}}\right)}
\end{equation}
 and
 \begin{equation} Y_UA_\zeta(-\alpha-\beta,\alpha-a;c-\alpha,\alpha)
= \frac{\zeta(1-a-\beta)\zeta(1+c)}{\zeta(1-\beta)\zeta(1+c-a)}
\prod_p \frac{\left(1-\frac{1}{p^{1+c}}\right)
\left(1-\frac{1}{p^{1+c-a}}
-\frac{1}{p^{1-\beta}}+\frac{1}{p^{1+c}}\right)}
{\left(1-\frac{1}{p^{1+ c-a}}\right)
\left(1-\frac{1}{p^{1-\beta}}\right)} .
\end{equation}
 If we differentiate the first expression with respect to $\beta$ and set $\beta=0$ we get
\begin{equation}
\frac{\zeta'}{\zeta}(1+a)-\frac{\zeta'}{\zeta}(1+c)-\sum_p\frac{\log p}{p^{2+a+c}}\left(1-\frac{1}{p^{1+c}}-\frac{1}{p^{1+a}}\right)^{-1}\left(
1-\frac{1}{p^{1+c}}\right)^{-1} .
\end{equation}
 The second expression gets multiplied by $(\frac{t}{2\pi})^{-a-\beta}$; then we differentiate with respect to $\beta$ and set $\beta=0$. Note that because of the factor $\zeta(1-\beta)$ in the denominator we only have to differentiate that term since it gives 0 when $\beta=0$. Thus,  the second term results in a contribution of
\begin{equation}
-\left(\frac{t}{2\pi}\right)^{-a}\frac{\zeta(1-a )\zeta(1+c)}{ \zeta(1+c-a)}
\prod_p \frac{\left(1-\frac{1}{p^{1+c}}\right)
\left(1-\frac{1}{p^{1+c-a}}
-\frac{1}{p }+\frac{1}{p^{1+c}}\right)}
{\left(1-\frac{1}{p^{1+ c-a}}\right)
\left(1-\frac{1}{p }\right)}  .
\end{equation}

Thus, we have
\begin{conjecture}Let $a$ and $c$ have non-negative real parts and satisfy $|a|,|c|\gg (\log T)^{-1}$. Then
\begin{align}
\frac{1}{T}\sum_{\gamma<T}\frac{\zeta(\rho+a)}{\zeta(\rho+c)}=&
\frac{1}{2\pi}\int_0^T\bigg(\log \frac{t}{2\pi}+
\frac{\zeta'}{\zeta}(1+a)-\frac{\zeta'}{\zeta}(1+c)\cr
&-\sum_p\frac{\log p}{p^{2+a+c}}\left(1-\frac{1}{p^{1+c}}-\frac{1}{p^{1+a}}\right)^{-1}\left(
1-\frac{1}{p^{1+c}}\right)^{-1}\cr
&\qquad -
\left(\frac{t}{2\pi}\right)^{-a}\frac{\zeta(1-a )\zeta(1+c)}{ \zeta(1+c-a)}
\prod_p \frac{\left(1-\frac{1}{p^{1+c}}\right)
\left(1-\frac{1}{p^{1+c-a}}
-\frac{1}{p }+\frac{1}{p^{1+c}}\right)}
{\left(1-\frac{1}{p^{1+ c-a}}\right)
\left(1-\frac{1}{p }\right)}\bigg)~dt \cr
&+O(T^{1/2+\epsilon}) .
\end{align}
\end{conjecture}

This agrees with Farmer's Conjecture (\ref{eqn:farmerconj}).

More applications of the ratios conjecture are given by Conrey and Snaith~\cite{CS}.


\end{document}